\documentclass[11pt,a4paper,reqno]{article}

\usepackage[
    a4paper,
    left=2.6cm,
    right=2.6cm,
    top=2.7cm,
    bottom=2.7cm
]{geometry}

\usepackage[T1]{fontenc}
\usepackage[utf8]{inputenc}
\usepackage{lmodern}
\usepackage{microtype}

\usepackage{amsmath}
\usepackage{amssymb}
\usepackage{amsthm}
\usepackage{mathrsfs}
\numberwithin{equation}{section}

\usepackage{graphicx}
\usepackage[usenames,dvipsnames]{xcolor}
\usepackage{float}
\usepackage{tikz}
\usetikzlibrary{decorations.pathmorphing}

\usepackage{authblk}
\let\affiliation\affil

\usepackage{setspace}
\usepackage{titlesec}
\titleformat{\section}
  {\normalfont\large\bfseries}
  {\thesection}{1em}{}
\titleformat{\subsection}
  {\normalfont\normalsize\bfseries}
  {\thesubsection}{1em}{}

\usepackage[font=small,labelfont=bf]{caption}

\usepackage[numbers,sort&compress]{natbib}
\usepackage[
    colorlinks=true,
    linkcolor=blue,
    citecolor=blue,
    urlcolor=blue
]{hyperref}

\DeclareMathOperator{\CT}{CT}
\DeclareMathOperator{\Res}{Res}

\newcommand{\eq}[1]{(\ref{#1})}
\newcommand{\nn}{\nonumber}

\newcommand{\beq}{\begin{equation}}
\newcommand{\eeq}{\end{equation}}
\newcommand{\beqq}{\begin{equation*}}
\newcommand{\eeqq}{\end{equation*}}
\newcommand\beqa{\begin{eqnarray}}
\newcommand\eeqa{\end{eqnarray}}
\newcommand\beqaa{\begin{eqnarray*}}
\newcommand\eeqaa{\end{eqnarray*}}
\newcommand\bea{\begin{array}}
\newcommand\eea{\end{array}}

\newtheorem{theorem}{Theorem}[section]
\newtheorem{definition}[theorem]{Definition}
\newtheorem{lemma}[theorem]{Lemma}
\theoremstyle{remark}

\title{Quantum Periods from TBA and Continued Fractions}

\author[1]{Alexander Chervov}
\author[2]{Ilia Gaiur}
\author[3,4,5]{Sergey Galkin}
\author[6]{Fedor Levkovich-Maslyuk}
\author[7, 8]{Volodya Rubtsov}
\author[8, 9]{Shamil Shakirov}

\affiliation[1]{Institute Imagine, Paris, France}

\affiliation[2]{Institut des Hautes \'Etudes Scientifiques\\
35 route de Chartres, 91440 Bures-sur-Yvette, France}

\affiliation[3]{Hetao Institute of Mathematics and Interdisciplinary Sciences (HIMIS), Shenzhen}

\affiliation[4]{Zhizhen School of Mathematics and Interdisciplinary Sciences (ZSIMS), CUHK, Hong Kong}

\affiliation[5]{on leave from PUC-Rio, Rio de Janeiro}

\affiliation[6]{Centre for Mathematical Science, City St George's, University of London, Northampton Square, EC1V 0HB, London, UK
 }

\affiliation[7]{LAREMA, UMR 6093 du CNRS, D\'epartement de Math\'ematiques,\\
Universit\'e d'Angers, 2, boulevard Lavoisier, 49045 Angers, France}

\affiliation[8]{Institute for Information Transmission Problems RAS, 127051 Moscow, Russia}

\affiliation[9]{ITMO University, St. Petersburg, Russia}

\begin{document}

\maketitle

\begin{abstract}
We consider quantum A-periods of local ${\mathbb P}^1 \times {\mathbb P}^1$, or equivalently the generating function of constant terms for powers of the Laurent polynomial $x + x^{-1} + y + y^{-1}$ where $x$ and $y$ are $q$-commuting variables. We derive a number of new representations for these rich combinatorial and algebraic objects. First, we find that the well-known TBA equation encodes not only the A-period itself but also 
%the coefficients of all diagonal terms generated by the noncommutative Laurent polynomials
higher diagonal modes of the corresponding noncommutative Laurent polynomials. Second, we obtain new formulas for the period in terms of continued fractions and related orthogonal polynomials, as well as a nontrivial explicit recursion allowing its efficient computation. Finally, we compute quantum corrections to the classical Picard-Fuchs equation, and observe that  holonomicity of the classical period does not generalise to $q$-holonomicity of its quantum counterpart. Applications to random walks on the Cayley graph of the Heisenberg group, and related new approaches to more efficient generative AI systems, are outlined.
\end{abstract}

\section{Introduction}

%\cite{AganagicEtAl2012} ,\cite{GRZZ25} %,\cite{HatsudaKatsuraTachikawa2016} - fundamental papers for %introduction

Quantum periods originate from semiclassical WKB expansions in quantum mechanics, and can be viewed as a quantisation of classical period integrals \cite{Dunham1932,Voros1983}. These are rich algebraic objects  that play a key role in a wide range of subjects within modern mathematical physics, including $\mathcal{N}=2$ supersymmetric gauge theories and AGT correspondence \cite{SeibergWitten1994,Nekrasov2003,NekrasovShatashvili2010,MironovMorozov2010SU2,MironovMorozov2010SUN,AldayGaiottoTachikawa2010}, quantum integrable systems and quantum spectral curves \cite{GorskyEtAl1995,Zenkevich2011,GukovSulkowski2012,Zenkevich2018}, refined topological string theory, mirror symmetry and the TS/ST correspondence \cite{AganagicEtAl2012,HuangKlemm2012,HuangEtAl2015,HatsudaEtAl2014,GrassiHatsudaMarino2016,GuEtAl2015,CodesidoGrassiMarino2017,Marino2018, KashaniPoor2016,KashaniPoorTroost2015,KashaniPoorTroost2014}, as well as resurgence, TBA equations and the ODE/IM duality \cite{GrassiGuMarino2020,ItoMarinoShu2019,DoreyDunningTateo2007}. At the same time, they incorporate deep mathematical structures related to combinatorics, algebraic and enumerative geometry, topological recursion, quantum cluster algebras, BPS wall-crossing and lattice spectral theory \cite{NekrasovOkounkov2006,EynardOrantin2007,BouchardEtAl2009,IwakiNakanishi2014,GaiottoMooreNeitzke2013,HatsudaKatsuraTachikawa2016}. 
%In this paper we uncover several new representations for %these important and widely used functions, motivated in %particular by a search for putative new structures such as a %direct quantisation of the classical Picard-Fuchs equation %for periods. Let us outline below the main ideas and %approaches we use.

Our algebraic approach is directly motivated by the work of
Gr\"afnitz, Ruddat, Zaslow and Zhou~\cite{GRZZ25}.
Their identification of the wavefunction definition of the
quantum A-period with the constant term of the logarithm
of the quantum mirror curve
\cite[Theorem~2.24]{GRZZ25}
provides the conceptual starting point for our study.

For the equal-coefficient specialization of local
$\mathbb{P}^{1}\times\mathbb{P}^{1}$, we prove that the formal
solution of the associated TBA functional equation generates
all diagonal coefficients of the quantum-torus resolvent,
with explicit normalization factors.
We develop this correspondence through continued fractions,
coefficient recursions and formal orthogonal polynomials.
We also investigate perturbative corrections to the classical
Picard--Fuchs equation and the possible existence of a scalar
linear $q$-difference equation in the modulus.
We begin by recalling the geometric and algebraic descriptions
of the periods under consideration.
\medskip 

\noindent
\textbf{Periods and quantization.} This paper is focused on  quantum versions of period integrals. We will stick to a very straightforward quantization procedure that mimics the standard Dirac quantization procedure for $(\mathbb{C}^*)^2$. Choices such as operator ordering and polarization enter
the construction and representation of the quantum Hamiltonian.
Boundary conditions enter the formulation of a spectral problem.
For a fixed quantum Hamiltonian and a chosen WKB branch,
period integrals arise in the formal continuation of the
WKB solution along cycles.

Concretely, for a state with energy $E$ the leading term of the WKB expansion $S_0(x,E)$ for the wavefunction can be formulated as an open-path type of integral of the Liouville form $\theta(E)$ restricted to the chosen level set of the classical system (i.e. in our case -- to a particular value of the energy). The standard WKB ansatz gives an expansion of the form
\begin{equation*}
    \Psi(x,E) \sim
    e^{\frac{1}{\hbar}\int\limits_{x_0}^x \theta(E)}
    \left(a_0(x,E)+\hbar a_1(x,E)+\hbar^2a_2(x,E)+\dots\right).
\end{equation*}
At leading exponential order in the WKB expansion,
formal continuation along a closed cycle
$\gamma\subset\{H=E\}$ is described by
\begin{equation}
    \Psi\Big\vert_{\gamma}
    \sim
    \Psi(x,E)
    e^{\frac{1}{\hbar}\oint\limits_\gamma \theta(E)}.
\end{equation}
The leading exponential continuation factor is therefore
determined by the period
\[
    \oint_{\gamma(E)} \theta(E).
\]
For a specified differential and a locally transported cycle,
with any required branch choices fixed, this period is a
locally defined function of the energy $E$.

The appearance of such periods in quantum spectral problems goes back to the semiclassical quantization conditions of WKB theory and their higher-order extensions \cite{Dunham1932}, and plays a central role in complex and exact WKB analysis \cite{Voros1983}. More recently, relations between quantum spectral problems and Bethe Ansatz equations have appeared in several forms, notably in the ODE/IM correspondence \cite{DoreyTateo1999,BazhanovLukyanovZamolodchikov2001} and in the quantization of classical integrable systems \cite{NekrasovShatashvili2010}. From this perspective, it is natural to ask how classical period integrals themselves should be quantized, and to what extent the resulting quantum periods can be described using methods of integrable systems and the Bethe Ansatz.

%Higher orders of the WKB expansion give formal quantum
%corrections to these classical periods.
%In the present paper, we study the quantum A-period through
%its constant-term description in the quantum torus
%$yx=qxy$. Our principal constructions use formal power
%series in the modulus, with coefficients depending on $q$.

Higher orders of the WKB expansion give formal quantum
corrections to classical periods.
For quantum mirror curves, a direct connection with the
constant-term approach is established by
Gr\"afnitz, Ruddat, Zaslow and Zhou~\cite{GRZZ25}. In their conventions, let $p(X,Y)=1-f(X,Y)$ denote the
quantum mirror curve, where $YX=qXY$, acting by
\[
    (X\Phi)(X)=X\Phi(X),
    \qquad
    (Y\Phi)(X)=\Phi(qX).
\]
A wavefunction satisfies $p\cdot\Psi=0$.
With the formal branch and normalization specified in
\cite{GRZZ25}, the wavefunction definition of the quantum
A-period reads
\[
    \underline a
    =
    \operatorname{CT}_{X}
    \log\frac{\Psi(qX)}{\Psi(X)}.
\]
Theorem~2.24 of \cite{GRZZ25} identifies this expression with
\[
    \underline a
    =
    \operatorname{CT}_{X,Y}\log p(X,Y)
    =
    -\sum_{n\geq1}\frac{1}{n}
       \operatorname{CT}_{X,Y}f(X,Y)^n.
\]
This identification directly motivates our algebraic
approach: it allows us to study quantum periods through
constant-term calculations in the quantum torus while
retaining their connection with the wavefunction description.

\medskip
\noindent
\textbf{Families of algebraic varieties in $(\mathbb{C}^*)^2$.} As we mentioned before, in this paper we study quantization of the period for families in $(\mathbb{C}^*)^2$. In particular, we study periods for the following family of algebraic curves
\begin{equation}\label{eq:Et}
    X_T: \quad x+x^{-1}+y+y^{-1} = 1/T.
\end{equation}
This is just a particular example of Laurent polynomials that appear in the classification of Fano varieties via mirror symmetry \cite{coates2014mirror,givental1999tutorial}. More precisely, this one is dual to $\mathbb{P}^1\times \mathbb{P}^1$, since the period of the holomorphic form $\gamma_t\in H^{1,0}(E_t)$  coincides with the quantum $Y$-function for the second Dubrovin connection associated with quantum cohomology $\text{QH}(\mathbb{P}^1\times \mathbb{P}^1)$. For more details, we refer to \cite{givental1999tutorial}. 

Although mainly we concentrate on methods not directly related to Fano or mirror symmetric aspects of the family $X_T$, initially our work was inspired by the technique commonly used in the mirror community. For families in $(\mathbb{C}^*)^n$ this technique allows one to write a period (at least locally) in a purely combinatorial way. In particular, for any Laurent polynomial $W(x_1,x_2,\dots x_n)\in \mathbb{Z}[x_1^{\pm 1},x_2^{\pm 1}, \dots x_n^{\pm 1}]$ the period for the family $W=1/t$ reads
\begin{equation*}
    \Pi(T):=\frac{1}{(2\pi i)^n}\oint\limits_{|x_i|=1} \frac{1}{1-TW}\bigwedge_{j=1}^n\frac{dx_j}{x_j}.
\end{equation*}
We give a more detailed explanation of this approach in Section \ref{sect:q=1}. 

\medskip
\noindent
\textbf{Combinatorial description.} Expanding this integral in the $t$ variable provides an extremely useful combinatorial description, namely
\begin{equation}\label{eq:CT}
    \Pi(T) = \sum\limits_{k=0}^{\infty} \CT_{x,y}[W^k] T^k,
\end{equation}
where $\CT$ stands for the constant term in corresponding variables. In the case of the $X_T$-family, the corresponding Laurent polynomial is $W=x+x^{-1}+y+y^{-1}$, and the only nonzero constant terms come from even powers so the period $\Pi(T)$ reads
\begin{equation}
    \Pi(T) = \sum\limits_{k=0}^{\infty} \binom{2k}{k}^2 T^{2k} = 1+ 4T^2 + 36T^4+\dots
\end{equation}
More generally, for a Laurent polynomial $f$ in commuting
variables, the classical period $\Pi_f(T)$ satisfies a
nonzero linear differential equation with polynomial
coefficients in $T$, and hence is $D$-finite.
This follows from the finite-dimensional Gauss--Manin
system governing the period; see Section~\ref{sect:q=1} and
\cite[Theorem~3.2]{coates2014mirror}.

%\begin{equation}
%    W=x+x^{-1}+y+y^{-1}
%\end{equation}
%Let us also mention that numbers appearing as coefficients in the %expansion are related to random walks on the square lattice with %uniform probability at each vertex. In particular, $[W^k]_0$ equals %the number of closed paths of length $k$ with fixed origin. This %algebraic description recovers a lot of stochastic properties of %such %random walks and may be also applied to more complicated %geometry of %the lattices. 

For the polynomial $W=x+x^{-1}+y+y^{-1}$ considered here,
the coefficients have a direct lattice-path interpretation.
The four monomials record the four nearest-neighbour steps
on $\mathbb Z^2$, and
$\operatorname{CT}_{x,y}W^k$ counts the walks of length $k$
that start and end at the origin.
For the simple symmetric random walk, each such walk has
probability $4^{-k}$, so the return probability is
\[
    p_k=4^{-k}\operatorname{CT}_{x,y}W^k.
\]
Every step changes the parity of the sum of the two
coordinates, which explains why the odd coefficients vanish.
This interpretation provides a combinatorial starting point
for the quantum deformation considered below.

\medskip
\noindent
\textbf{Noncommutative variables and quantum periods.} 
While the formula \eqref{eq:CT} has a deep geometric origin, for us its main feature is the possibility to recover the period function without explicitly referring to the geometry of the family $X_T$. This allows us to pass to the non-commutative fiber straightforwardly sending $x$ and $y$ to any noncommutative algebra. The same holds for an arbitrary Laurent polynomial. The case of constant term elements of Laurent polynomials taking values in free algebras was considered in \cite{kontsevich2011noncommutative}. Remarkably, the author shows that such generating functions are in fact algebraic.

The constant-term description of the classical period suggests
an algebraic approach to quantization. Specifically, since $X_{T}$ is a family of curves in $(\mathbb{C}^*)^2$, the natural procedure is to quantize corresponding volume form $\omega_{(\mathbb{C}^*)^2} = \frac{dx\wedge dy}{xy}$.  Following the standard deformation quantization procedure \cite{dirac1925fundamental, kontsevich2003deformation}, one gets a noncommutative algebra known as quantum torus
\[
\mathcal K
=
\mathbb C[q,q^{-1}]
\langle x^{\pm1},y^{\pm1}\rangle
\big/(yx-qxy),
\]
where $q$ is a deformation parameter that is central. In $\mathcal{K}$ every element has a unique
finite expansion
\[
U=\sum_{a,b\in\mathbb Z}c_{a,b}(q)x^a y^b,
\]
and one can introduce trace as
\[
\operatorname{CT}_{x,y}U:=c_{0,0}(q).
\]
Indeed one checks that it is a trace functional, i.e. $\operatorname{CT}_{x,y}(UV) = \operatorname{CT}_{x,y}(VU)$. It obviously extends to the ring of formal power series $\mathcal K[[T]]$. The algebra has following representation via operators
\[
(x\phi)(u)=u\phi(u),
\qquad
(y\phi)(u)=\phi(qu).
\]
This allows to send the functions on $(\mathbb{C}^*)^2$ to explicit operators. For
$W$ the corresponding spectral problem reads as
\begin{equation}\label{eq:qdiff}
\Psi(qu)+\Psi(q^{-1}u)
+\bigl(u+u^{-1}\bigr)\Psi(u)
=E\Psi(u),
\qquad E=T^{-1}.
\end{equation}
The generating function for constant terms reads as
\begin{equation*}
    G(T) = \sum\limits_{k=0}^{\infty} \CT_{x,y}[W^k] T^k \in \mathbb{Z}[q^{\pm 1}][[T]].
\end{equation*}
%While $G(T)$ is not usual for theoretical physics literature, its %Euler antiderivative is known as the quantum A-period of $W$. It is %defined as 
%\begin{equation}
    %F(T) = -\sum\limits_{j=1}^{\infty}\CT_{x,y}[W^k] \frac{T^k}{k} = %\frac{1}{T}\int G(T) dT - 1.
%\end{equation}
The corresponding quantum A-period, in the constant-term
normalization used here, is
\begin{equation}
\begin{aligned}
F(T)
:=\operatorname{CT}_{x,y}\log(1-TW)
=-\sum_{k\ge1}
  \frac{\operatorname{CT}_{x,y}(W^k)}{k}\,T^k.
\end{aligned}
\end{equation}
A-period is related to the constant term generating function as follows
\[
G(T)=1-T\frac{dF}{dT}.
\]
Quantum spectral problem \eqref{eq:qdiff} has a direct interpretation via the Harper--Hofstadter model \cite{avila2017spectral,jitomirskaya2022one}.
For $q=e^{2\pi i\alpha}$, with $\alpha\in\mathbb R$, the unitary representation on $\ell^2(\mathbb Z)$ defined by
\[
(xu)_n=e^{i(\theta+2\pi\alpha n)}u_n,
\qquad (yu)_n=u_{n+1}
\]
satisfies $yx=qxy$ and realizes $W$ as the isotropic Harper operator
\begin{equation}\label{eq:harper-realization}
(H_{\alpha,\theta}u)_n
=u_{n+1}+u_{n-1}+2\cos(\theta+2\pi\alpha n)u_n.
\end{equation}
The relation between magnetic translations, quantum groups and Bethe-Ansatz equations for the mid-band problem was developed by Wiegmann and Zabrodin~\cite{WiegmannZabrodin1994}.
In this representation, taking the matrix element at $\delta_0$ and then averaging over $\theta$ recovers $\operatorname{CT}_{x,y}$: for $x^a y^b$ the matrix element is $e^{ia\theta}\delta_{b,0}$.
Let $d\nu_\alpha$ denote the normalized density-of-states measure, obtained by this phase average of the spectral measure at $\delta_0$.
Since $\|H_{\alpha,\theta}\|\leq4$, for $|T|<1/4$ we obtain
\begin{equation}\label{eq:harper-periods}
\begin{aligned}
G(T)
&=\frac1{2\pi}\int_0^{2\pi}
\langle\delta_0,(1-TH_{\alpha,\theta})^{-1}\delta_0\rangle\,d\theta
=\int_{\mathbb R}\frac{d\nu_\alpha(\lambda)}{1-T\lambda},\\
F(T)&=\int_{\mathbb R}\log(1-T\lambda)\,d\nu_\alpha(\lambda),
\end{aligned}
\end{equation}
where the logarithm is the branch vanishing at $T=0$ and $q=e^{2\pi i\alpha}$ is understood in $G$ and $F$.
Thus our moment and continued-fraction formulas compute the large-energy expansion of the averaged Green function $E^{-1}G(E^{-1})$.
The relation between quantum A-periods and the Hofstadter density of states was established at rational magnetic flux in~\cite[Section~3.3, equation~(80)]{HatsudaKatsuraTachikawa2016}; our emphasis is on its algebraic realization through TBA and on the correspondence for all diagonal modes.

%The function appears prominently in modern mathematical physics in %several contexts, including WKB quantization for Hamiltonians with %exponential terms, and quantum geometry of Calabi-Yau's. For the %quantum mechanical Hamiltonian \eqref{eq:qdiff}, the parameter $q$ %naturally appears as the exponentiated Planck constant. The WKB %wavefunction of the Hamiltonian can be constructed as a small %parameter expansion, and then the quantum period $F(T)$ appears as a %contour integral (constant term) of either the derivative or the $q$-%derivative of the logarithm of the wavefunction \cite{B24,GRZZ25}. %The equivalence of this definition and the simple algebraic %definition via the constant term that we use is proven in %\cite{GRZZ25} (see Th. 2.24).
The relation between the constant-term construction and the
wavefunction definition of quantum periods was recalled above.
Their equivalence, with the formal branch and normalization
specified there, is established in
\cite[Theorem~2.24]{GRZZ25}.
In our algebraic calculations, $q$ is an indeterminate;
its interpretation as an exponentiated Planck constant arises
in the wavefunction realization.

Going back to mirror symmetry, one thinks of \eqref{eq:qdiff} as of the quantum curve for local $\mathbb{P}^1\times \mathbb{P}^1$. This is closely related to the WKB context, the difference operators being the exponentiated momentum, or shift. Spectral properties of the quantum curve operator are important in the context of the Topological String/Spectral Theory correspondence \cite{CodesidoGrassiMarino2017,GrassiHatsudaMarino2016,Marino2018}. The function $F(T)$ appears as a quantization of the period integral for the mirror curve, the quantum $A$-period of the Calabi-Yau \cite{AganagicEtAl2012}. For an explicit example see equation (5.30) in \cite{GrassiHatsudaMarino2016}. For comparison, setting $T=e^{-\mu}$, the terms displayed
in equation~(5.30) of \cite{GrassiHatsudaMarino2016} agree with $2F(T)$,
reflecting a different normalization of the quantum period.

Furthermore, it is clear that while this paper works with a specific choice of Laurent polynomial, much of the above discussion can be generalized to more general $W$. For example, the choice $x + y + x^{-1} y^{-1}$ is relevant to the quantum geometry of the Calabi-Yau local ${\mathbb P}^2$, and is thoroughly studied. Table 1 of \cite{GrassiHatsudaMarino2016} gives more examples. Another generalization which is not difficult to implement is introducing non-trivial coefficients in the polynomial $W$ -- for example, to consider $m( x + x^{-1} ) + y + y^{-1}$ instead of $x + x^{-1} + y + y^{-1}$.

\medskip

\noindent 
%\textbf{TBA equations.} 

% Dirac quantization procedure sends classical Hamiltonian $H(\mathbf{q},\mathbf{p})$ written in the Darboux variables $\omega(\mathbf{q},\mathbf{p}) = \sum dp_i\wedge d q_i$ to the operator $\hat{H}$ by sending
% \begin{equation*}
% p_i\,\rightarrow\,\hat{p}_i=\hbar\frac{\partial}{\partial x_i},\quad q_i\,\rightarrow\,\hat{q}_i=x_i\cdot,\quad [\hat{p}_i,\hat{q}_i]=\hbar.
% \end{equation*}

\medskip

\paragraph{TBA equation and main results.}
Relations between quantum A-periods and functional equations
associated with the Thermodynamic Bethe Ansatz have been studied
in \cite{HatsudaEtAl2014,HuangTBA,DuHuang}.
For the equal-coefficient specialization considered here,
the relevant TBA-like equation is
\begin{equation}
\label{eq:intro-tba}
\begin{aligned}
\eta(T,z)^2
&=1+T^2 S(z)
  \bigl(\eta(T,qz)+\eta(T,z)\bigr)\\
&\qquad{}\times
  \bigl(\eta(T,q^{-1}z)+\eta(T,z)\bigr),
\end{aligned}
\end{equation}
where $S(z)=2+z+z^{-1}$.
This is the $m=0$ specialization of
\cite[equation~(2.8)]{HuangTBA}, after renaming the auxiliary
variable and replacing its modulus by $T^2$.
We study this equation formally in $T$, with $q$ an indeterminate.
Its shifts act on the auxiliary variable $z$. The equation has a unique solution
\[
\eta(T,z)\in
1+T^2\mathbb C(q)[z,z^{-1}][[T^2]].
\]
Our main result identifies every Laurent coefficient of this
solution with a diagonal coefficient of the quantum-torus
resolvent. For $W=x+x^{-1}+y+y^{-1}$ with $yx=qxy$, we prove
\begin{equation}
\label{eq:intro-main-result}
[z^m]\eta(T,z)
=
\frac{2q^{-m(m-1)/2}}{1+q^m}\,
[x^m y^m](1-TW)^{-1},
\qquad m\in\mathbb Z.
\end{equation}
Here coefficient extraction on the right is taken in the
ordered monomial basis $x^a y^b$.
This identity is proved in Section~\ref{sect:TBA}. Its zero-mode specialization is
\[
\operatorname{CT}_z\eta(T,z)
=
\operatorname{CT}_{x,y}(1-TW)^{-1}
=
G(T)
=
1-T\frac{dF}{dT}.
\]
Together with $F(0)=0$, this recovers the quantum A-period
by formal integration. The nonzero modes provide the
corresponding generating functions for all remaining
diagonal coefficients.
The proof combines identities in the quantum torus with
the linearized TBA equation and the normalization
$\eta(T,-1)=1$.

We develop this correspondence through continued-fraction
representations, lattice-path decompositions and explicit
coefficient recursions. In particular, the recursion
of Section~\ref{sect:methodII} computes $\eta(T,z)$ through order $T^{2N}$
using $O(N^2)$ arithmetic operations over $\mathbb C(q)$.
The resulting three-term recurrences also give rise
to formal orthogonal polynomials.

Finally, we compute perturbative corrections to the
classical Picard--Fuchs equation and investigate scalar
linear $q$-difference equations in the modulus $\tau=T^2$.
Our search finds no such equation of order at most eight
with polynomial coefficients of degree at most eighteen.
The existence of an equation beyond these bounds
remains open.

\pagebreak

\paragraph{Cayley graphs, random walks and connections to AI}
Consider the discrete Heisenberg group
\[
\mathsf H=\langle x,y,q\mid yx=qxy,\, xq=qx,\, yq=qy \rangle
\]
with the symmetric generating set $\{x,x^{-1},y,y^{-1}\}$.
In the group algebra, right multiplication by $\widetilde W=x+x^{-1}+y+y^{-1}$ gives the adjacency operator of its Cayley graph.
The quantity $N_n(r,s,k)=[x^ry^sq^k]\widetilde W^n$,  counts length-$n$ walks from the identity to $x^ry^sq^k$, and their transition probability is $4^{-n}N_n(r,s,k)$. In particular, $\CT_{x,y}W^n$ counts walks whose projection to $\mathbb Z^2$ is closed, with $q$ recording the central displacement, while the coefficient of $q^0$ counts returns to the identity in $\mathsf H$.

These observations are also related to some AI approaches, in the following way. For Cayley graphs each edge is labelled by the generator, and hence paths on a graph define a word in the generators (see \cite{CayleyPyCube2025,CayleyPyRL2026,CayleyPyGrowth2025,CayleyPyHolography2026,DouglasFraserTaliente2025,HeYau2022,McAleerEtAl2019DeepCube,AgostinelliEtAl2019DeepCubeA,Takano2023EfficientCube,Johnson2018RubiksGuidance,BarketEtAl2026CayleyGraphs,WeiEtAl2026Dihedral}
% McAleerEtAl2019DeepCube,AgostinelliEtAl2019DeepCubeA,Takano2023EfficientCube,Johnson2018RubiksGuidance,BarketEtAl2026CayleyGraphs,WeiEtAl2026Dihedral,
for recent machine learning-based explorations of Cayley graphs of immense size and related objects). Considering various admissible classes of paths one obtains various sets of words, i.e. languages.
The random walks defined by $W^n$ are one example. 
Assigning more non-uniform probabilities to edges one obtains other languages.
The task of generative AI is to describe the probabilities of generating different words starting from some initial word. 
That means, essentially, to describe analogues of $W^n$. For Cayley graphs of the commutative groups $W^n$ satisfy recursive Picard-Fuchs equations which allow their efficient computation, and one motivation for the present paper is to explore how similar effective computations can be extended beyond the commutative groups,
and thus potentially generalized to practical AI tasks leading to more effective computational frameworks. The algorithm eq. (\ref{Algorithm2}) that we found appears to be the best substitute of Picard-Fuchs equations beyond the commutative case. 

\medskip
\noindent
\paragraph{Organization of the paper.}
Section~\ref{sect:q=1} reviews the classical limit, its elliptic-curve
interpretation and the associated Picard--Fuchs equation.
Section~\ref{sect:A-period} introduces the quantum-torus formulation of the
A-period and the generating function for diagonal coefficients.
Section~\ref{sect:TBA} proves their correspondence with the formal TBA
solution. Section~\ref{sect:methodI} develops continued-fraction representations,
while Section~\ref{sec:barrier} presents an alternative combinatorial proof
of the full correspondence using lattice paths and a diagonal
barrier, together with an independent proof of the identity for the zero mode. Section~\ref{sect:methodII} derives explicit coefficient recursions
and an efficient computational algorithm, and Section~\ref{sect:methodIII}
studies the associated formal orthogonal polynomials.
Section~\ref{sect:PF} examines perturbative corrections to the classical
Picard--Fuchs equation and searches for scalar linear
$q$-difference equations in the modulus.
Section~\ref{sect:conclusions} discusses the conclusions and further directions.
Appendix~\ref{app:ode} records the perturbative differential equation
to higher order.

\paragraph{Statement of the generative AI usage.}  During the preparation of this work the authors used various services, models and software packages (from ChatGPT, Claude, Google DeepMind and Every) to transcribe voice recordings and blackboard photos into text, check the language and text structure, literature search and review, as well as an assistant tool during the proofs and computational routines. After using this tool/service, the authors reviewed and edited the content as needed and take full responsibility for the content of the published article. We thank Google DeepMind for providing the access to some of these tools.

\pagebreak

\section{The classical limit  $q=1$}\label{sect:q=1}

This section contains a review of the commutative case $xy=yx$. The main aim is to explain why this case has a geometric origin and as a consequence gives holonomic functions. Before we go into details, we remind the reader the definitions of a holonomic function and discuss periods of algebraic varieties.

\begin{definition} 

A series is called \emph{$D$-finite} (holonomic) over a field
$L$ of characteristic zero if it satisfies a nonzero homogeneous linear
differential equation with coefficients in $L[T]$. Equivalently, its
coefficients satisfy a linear recurrence with polynomial coefficients
in the index. This refers to ordinary differentiation, not to a
$q$-difference equation.

\end{definition}

\noindent
By ``geometric origin'' we mean that the function is represented by
integrating an algebraic differential form over a cycle in a family of
algebraic varieties. The resulting finite-dimensional cohomology carries
the Gauss--Manin connection; eliminating its other components gives a
Picard--Fuchs equation for the chosen period. Thus an ordinary algebraic
period must be $D$-finite. The converse is not being asserted. Over
$\overline{\mathbb Q}$ this is the setting of arithmetic geometry and
motives; see Y. Andr\'e~\cite{Andre89}. We use \emph{motivic origin}
only in this concrete sense of a geometric realization. In the rest of the paper we explicitly show how constant term generating series in the commutative case leads to the elliptic curves family period.

From the  naive semi-classical point of view, the problem drops to the consideration of level curves for the classical Hamiltonian, i.e. a family of algebraic curves in $\mathbb{C}^*\times \mathbb{C}^*$ given by
\begin{equation}
    x+x^{-1}+y+y^{-1} = E,\quad E\in \mathbb{C}.
\end{equation}
Extending $E$ to $\mathbb{CP}^1$ and choosing the north pole chart $E=1/T$, one ends up with the so-called Landau-Ginzburg model for the Laurent polynomial $W(x,y)$
$$
W(x,y)=x+x^{-1}+y+y^{-1}=1/T.
$$
In order to recover the classical limit of the quantum A-period, we consider the same generating function for the constant terms, which as it was mentioned before coincides with the elliptic integral $K(4t)$ up to a constant. This observation trivially follows from the fact that the generating function can be written as a contour integral of the form
\begin{equation}\label{eq:LG}
   \frac{1}{(2\pi i)^2}\oint\limits_{|x|=|y|=1} \frac{dxdy}{xy}\frac{1}{1-TW(x,y)} := \Pi(T),
\end{equation}
where $|T|<1/4$. Put \(S(x)=x+x^{-1}\) and
\[
G_t(x,y)=y(1-TW(x,y))
       =-Ty^2+(1-TS(x))y-T.
\]
For \(t\) sufficiently close to zero, the \(y\)-contour encloses the
root
\[
y_-(x,T)=
\frac{1-TS(x)-\sqrt{(1-TS(x))^2-4T^2}}{2T},
\]
where the square root tends to \(1\) as \(T\to0\). Taking the residue
at this root yields
\[
\Pi(T)=\frac{1}{2\pi i}\oint_{|x|=\varepsilon}
\frac{dx/x}{\sqrt{(1-T(x+x^{-1}))^2-4T^2}}.
\]
Equivalently, set
\[
v^2=\bigl(x-T(x^2+1)\bigr)^2-4T^2x^2
   =\bigl[x-T(x+1)^2\bigr]\bigl[x-T(x-1)^2\bigr].
\]
With the branch \(v\to x\) as \(T\to0\), the period takes the
form of elliptic integral
\[
\Pi(T)=\frac{1}{2\pi i}\oint_{\gamma_T}\frac{dx}{v}.
\]
For generic \(T\), this quartic has four distinct roots. Its smooth
projective completion is therefore a genus-one curve, and \(dx/v\)
is its holomorphic differential. Finally, expansion of the original constant-term integral gives
\[
\Pi(T)=\sum_{n\geq0}\binom{2n}{n}^{\!2}T^{2n}
={}_2F_1\!\left(\frac12,\frac12;1;16T^2\right)
=\frac{2}{\pi}K(4T),
\]
where \(K(k)=\int_0^{\pi/2}(1-k^2\sin^2\theta)^{-1/2}\,d\theta\)
uses the modulus convention.

Let us now give fully geometric description of the family. For generic $T$, the affine curve $X_T$ is smooth.
Let $E_T$ denote its smooth projective compactification,
which has genus one.
Write $p_1(T),\ldots,p_4(T)$ for the four distinct points
of $E_T\setminus X_T$, and define the reduced effective
boundary divisor
\[
    D_T:=p_1(T)+p_2(T)+p_3(T)+p_4(T),
    \qquad \deg D_T=4.
\]
Thus
\[
    X_T=E_T\setminus\operatorname{Supp}(D_T).
\]
Two projective realizations of $E_T$ are described below. The elliptic interpretation can be made explicit in two
different projective realizations.
Multiplying the defining equation of $X_T$ by $xy$ gives
\[
    F_T(x,y)
    :=x^2y+xy^2+x+y-T^{-1}xy=0.
\]
Its homogenization defines a plane cubic in $\mathbb P^2$:
\[
    u^2v+uv^2+uw^2+vw^2-T^{-1}uvw=0,
    \qquad x=u/w,\quad y=v/w.
\]
Alternatively, compactifying the two coordinates separately
gives a curve in $\mathbb P^1\times\mathbb P^1$.
With $x=u_1/u_0$ and $y=v_1/v_0$, its equation is
\[
\begin{aligned}
    &(u_1^2+u_0^2)v_0v_1
      +u_0u_1(v_1^2+v_0^2)\\
    &\qquad
      -T^{-1}u_0u_1v_0v_1=0.
\end{aligned}
\]
This equation has bidegree $(2,2)$.
In this model, the points supporting $D_T$ can be labelled as
\[
\begin{aligned}
    p_1(T)&=(0,0),       & p_2(T)&=(0,\infty),\\
    p_3(T)&=(\infty,0),  & p_4(T)&=(\infty,\infty).
\end{aligned}
\]
Their ambient coordinates are independent of $T$,
and each point occurs in $D_T$ with multiplicity one. Over
\[
    B=\mathbb P^1_T
      \setminus\{0,1/4,-1/4,\infty\},
\]
both projective models are smooth curves of genus one.
Their common open subset $X_T$ determines a birational
identification, which extends to an isomorphism of the
smooth projective curves.
Thus they are two projective realizations of the same
compactification $E_T$.
Choosing a boundary point as the origin makes $E_T$
an elliptic curve.

Let us note that the residue construction is not specific to the Laurent polynomial above. Let
\[
W\in\mathbb{Q}[x_1^{\pm1},\ldots,x_n^{\pm1}]
\]
be a nonconstant Laurent polynomial, such that corresponding family of algebraic surfaces
\[
X_T=\{(x_1,\ldots,x_n)\in(\mathbb{C}^*)^n: W(x)=1/T\},
\]
is smooth for generic value of $T$, i.e. $\partial/\partial x_k W = 0$ only for specific values of parameter $T$. The form
\[
\Omega_T=
\frac{1}{1-TW(x)}
\bigwedge_{j=1}^n\frac{dx_j}{x_j}
\]
has a simple pole along \(X_T\). Its Poincar\'e residue is a
holomorphic \((n-1)\)-form \(\omega_t\) on \(X_T\). For example, on a
patch where \(\partial W/\partial x_n\ne0\), using the convention
\(\Omega_T=\omega_T\wedge d(1-TW)/(1-TW)\), it is given by
\[
\omega_T=
-\left.
\frac{1}{T x_n\,\partial W/\partial x_n}
\bigwedge_{j=1}^{n-1}\frac{dx_j}{x_j}
\right|_{X_T}.
\]
Whenever the integration cycle can be expressed as a tube around a
cycle \(\gamma_T\subset X_T\), the residue theorem turns the
\(n\)-dimensional integral of \(\Omega_T\) into a period
\[
\frac{1}{(2\pi i)^{n-1}}\int_{\gamma_T}\omega_T.
\]
In two variables, \(X_T\) is a curve and \(\omega_T\) is a
holomorphic one-form, as in the preceding example. 

As we discussed before, the geometric interpretation also implies that constant term generating functions are holonomic. As \(T\) varies, the residue classes
\([\omega_T]\) lie in a finite-dimensional algebraic de Rham
cohomology space equipped with the Gauss--Manin connection.
Consequently, some finite collection of their successive
\(T\)-derivatives is linearly dependent over \(\mathbb{Q}(T)\).
Integrating this relation over a locally constant cycle gives
\[
\sum_{k=0}^{r} a_k(T)\,\frac{d^k}{dT^k}\Pi_W(T)=0,
\qquad a_k(T)\in\mathbb{Q}(T),
\]
where the \(a_k\) are not all zero. After clearing denominators, this is
a linear differential equation with polynomial coefficients; thus
\(\Pi_W(T)\) is holonomic \cite{BostanStraubYurkevich23}.

\medskip
\noindent 
\textbf{Relation to Bessel function.} The potential $W(x,y)$ that we consider closely related to the products of Bessel functions
\begin{equation}
    I_0(2s) = \oint e^{s (x+x^{-1})} \frac{dx}{x}.
\end{equation}
Expanding integrand into Taylor series, one relates the function with constant terms of the exponential argument. Taking its square one gets
\begin{equation}
I_0(2s)^2 = \oint\limits_{|x|=|y|=1} e^{s (x+x^{-1}+y+y^{-1})} \frac{dx}{x}\frac{dy}{y}
=\sum_{n\geq0}\binom{2n}{n}^{\!2}\frac{s^{2n}}{(2n)!}.
\label{BesselIdentity}
\end{equation}
This is the classical Bessel product mechanism; higher analogues are
studied in \cite{HigherBessel26}. Using $q=1$ for quantum periods introduced before
\begin{align}
F(T)\big|_{q=1}&=-\frac12\sum_{n\geq1}\binom{2n}{n}^{\!2}\frac{T^{2n}}n
=-2T^2\,{}_4F_3\!\left(\begin{matrix}1,1,\frac32,\frac32\\2,2,2\end{matrix};16T^2\right),\\
G(T)\big|_{q=1}&={}_2F_1\!\left(\frac12,\frac12;1;16T^2\right)
=\frac2\pi K(4T),
\end{align}
one relates these classical limits to Laplace transform of the squares of Bessel functions
\begin{align}
G(T)\big|_{q=1}&=\int_0^\infty e^{-r}I_0(2Tr)^2\,dr,\\
F(T)\big|_{q=1}&=-\int_0^\infty\frac{e^{-r}}r
                  \bigl(I_0(2Tr)^2-1\bigr)\,dr.
\end{align}
These formulas are specific case of Bessel identities appeared in the study of fullerenes and corresponding random walks on graphs \cite{bille2025random}.

\section{Quantum A-Period}\label{sect:A-period}

In this section we introduce all required mathematical structures and notions in order to define the quantum A-period and it's generalization to higher correlators.

Fix a non-zero complex parameter $q \in {\mathbb C}^{\star}$ and consider the  algebra $\mathcal{K}=\mathbb{C}\langle x^{\pm 1},y^{\pm 1}\rangle/\langle yx=qxy\rangle$, also known as quantum plane \cite{manin1987some}. One may choose a basis given by $e_{n,m}=x^ny^m$. Indeed, every monomial in $x,y,x^{-1},y^{-1}$ can be reduced up to some power of $q$
to one of $e_{n,m}$, using the relation $yx = qxy$ and its corollaries $$y^{-1}x=q^{-1}xy^{-1}, \quad yx^{-1} = q^{-1} x^{-1}y, \quad y^{-1}x^{-1}  = q x^{-1}y^{-1}.$$
The algebra $\mathcal{K}$ acts on itself by multiplication, and in the basis $e_{n,m}$ the action reads
\begin{equation*}
    x e_{n,m} = e_{n+1,m},\quad y e_{n,m} = q^n e_{n,m+1}.
\end{equation*}
Moreover, one may define a scalar product for these basis elements:
\begin{equation}
\big< e_{n,m} , e_{k,l} \big> = \delta_{n,k} \delta_{m,l}.
\end{equation}
The scalar product defines projectors (additive characters) that send an element of the algebra $\mathcal{K}$ to $\mathbb{C}[q^{\pm 1}]$ by taking  the coefficient in front of the corresponding monomial. More precisely, one has
$$
[L]_{x^iy^j} = \langle e_{i,j}, L\rangle,
\quad L\in \mathcal{K}.
$$
Clearly, the Laurent polynomial $W$ (as well as its powers) belongs to the algebra $\mathcal{K}$. Using projectors one can write the quantum A-period for $W$ in the form
\begin{align}
F(T) = \left< e_{0,0}, \ln(1 - T W) \right> := - \sum\limits_{n = 1}^{\infty} \left< e_{0,0}, W^n \right> \dfrac{T^n}{n}
\label{QuantumAPeriod}
\end{align}
together with the corresponding constant term generating function:
\begin{align}
G(T) = \left< e_{0,0}, \dfrac{1}{1 - T W} \right> := \sum\limits_{n = 0}^{\infty} \left< e_{0,0}, W^n \right> T^n = 1 - T \dfrac{d F(T)}{dT}
\label{QuantumAPeriodG}
\end{align}
We also introduce the generating function for higher matrix elements that generalizes $G(T)$ as follows:
\begin{align}
\eta(T,z) = \sum\limits_{n = 0}^{\infty} \sum\limits_{m = -\infty}^{\infty} \dfrac{ 2 q^{-m(m-1)/2} }{1 + q^{m} } \ \left< e_{m,m}, W^n \right> T^n z^m
\label{GenFuncHigher}
\end{align}
This specific form of coefficients is motivated by the connection with TBA we discuss below. As one can see, the function $G(T)$ can be obtained as the constant term in $z$ of $\eta(T,z)$. 

%Relation to graphs and AI.

% Finally, it is worth mentioning that the commutative case $q=1$ is related to the well-developed theory of classical periods. We have by substitution

% \begin{align}
% F(T)\Big|_{q=1} = \frac{-1}{2} \left( 4 T^2 + 18 T^4 + \dfrac{400}{3} T^6 + O( T^8 ) \right)
% \label{QuantumAPeriodSeriesCl1}
% \end{align}
% \smallskip\\
% and in fact this can be computed in closed form:

% \begin{align}
% F(T)\Big|_{q=1} = \frac{-1}{2} \sum\limits_{n = 1}^{\infty} \dfrac{T^{2n}}{n} \left(\dfrac{(2n)!}{n!^2} \right)^2 = 
% -2T^2\,\,{}_4F_3\left(
% \begin{matrix}
% 1 & 1 & \frac32 & \frac32\\
% 2 & 2 & 2
% \end{matrix}
% \,\middle|\,16T^2
% \right).
% % -2 T^2 \ _{4}F_{3}\left(\begin{array}{llllll} 1 & 1 & 3/2 & 3/2 \\ 2 & 2 & 2 \end{array} \big| 16 T^2\right)
% \label{QuantumAPeriodSeriesCl2}
% \end{align}
% \smallskip\\
% because the constant term of a power $P^{2n}$ of commutative polynomial is easily shown to be $((2n)!/n!^2)^2$. Perhaps more illuminating is the expression for the other generating function,

% \begin{align}
% \left< e_{0,0}, \dfrac{1}{1 - T P} \right> = \sum\limits_{n = 0}^{\infty} \left< e_{0,0}, P^n \right> T^n = \sum\limits_{n = 0}^{\infty} \left(\dfrac{(2n)!}{n!^2} \right)^2 T^{2n} = \dfrac{2}{\pi} K( 16 T^2 )
% \label{QuantumAPeriodSeriesCl3}
% \end{align}
% \smallskip\\
% directly expressed through the complete elliptic $K$ integral. In the next subsection we discuss the $q=1$ case in more detail.

% Kontsevich's ideas??? Commutative case is always holonomic??? Picard-Fuchs??? 

\section{Solving Thermodynamic Bethe Ansatz equation}\label{sect:TBA}

In this section we prove that the generating function of higher matrix elements \eqref{GenFuncHigher}
satisfies quadratic TBA equation. To the best of our knowledge, such relation was known only for A-period functions. First let's us introduce quadratic TBA equation and discuss uniqness of it's solution.
For fixed \(q\in\mathbb C^\times\), set
\[
\mathcal R_z:=\mathbb C[z,z^{-1}],
\qquad
S(z):=z+2+z^{-1}.
\]

\begin{theorem}[Formal TBA solution]\label{thm:TBA}
The equation
\begin{equation}
1-\eta(T,z)^2+T^2\bigl(\eta(T,qz)+\eta(T,z)\bigr)
\bigl(\eta(T,q^{-1}z)+\eta(T,z)\bigr)S(z)=0
\label{EquationTBA1}
\end{equation}
has a unique solution $\eta\in1+T^2\mathcal R_z[[T^2]]$.
Its coefficient of $T^{2n}$ has Laurent support in $[-n,n]$
(that is, only the powers $z^{-n},\ldots,z^n$ occur), and
$\eta(T,z)=\eta(T,z^{-1})$.
\end{theorem}
\begin{proof}
Write $\eta=1+\sum_{n\geq1}\eta_nT^{2n}$. At order $T^{2n}$ the only
new coefficient is $-2\eta_n$; the other terms involve earlier
coefficients. This gives existence, uniqueness and the support bound
inductively. The equation is invariant under $z\mapsto z^{-1}$,
which proves symmetry by uniqueness.
\end{proof}
\noindent
Equation \eqref{EquationTBA1} is ${\mathfrak m}=0$ specialization of equation appeared in 
\cite{HuangTBA} (see eq.~(2.8)), with its modulus replaced by $T^2$.
Its origin in thermodynamic Bethe ansatz (TBA) integral equations and
related functional relations is discussed in \cite{HatsudaEtAl2014,HuangTBA,DuHuang}.
Here only the formal difference equation is required.

Now we propose a useful technical result which we shall use later.

\begin{lemma}[The quadratic residual]
\label{lem:tba-residual}
Let $q$ be an indeterminate, set
\[
    S(z)=2+z+z^{-1},
\]
and suppose that
$\eta(T,z)\in\mathbb{C}(q)[z,z^{-1}][[T]]$
satisfies the linear identity
\[
\begin{aligned}
\eta(T,z)-\eta(T,qz)
=T^2\Bigl[
&S(z)\bigl(\eta(T,q^{-1}z)+\eta(T,z)\bigr)\\
&-S(qz)\bigl(\eta(T,qz)+\eta(T,q^2z)\bigr)
\Bigr].
\end{aligned}
\]
Define the quadratic residual by
\[
\begin{aligned}
\mathcal E(T,z)
={}&1-\eta(T,z)^2\\
&+T^2S(z)
\bigl(\eta(T,z)+\eta(T,qz)\bigr)
\bigl(\eta(T,z)+\eta(T,q^{-1}z)\bigr).
\end{aligned}
\]
Then
\[
    \mathcal E(T,qz)=\mathcal E(T,z),
\]
and consequently $\mathcal E(T,z)$ is independent of $z$.
If, in addition, $\eta(T,-1)=1$, then
$\mathcal E(T,z)=0$.
\end{lemma}

\begin{proof}
For brevity, write
\[
\begin{aligned}
a&=\eta(T,q^{-1}z), & b&=\eta(T,z),\\
c&=\eta(T,qz),     & d&=\eta(T,q^2z),
\end{aligned}
\qquad
s=S(z),\quad s_+=S(qz).
\]
The assumed linear identity becomes
\[
    b-c=T^2\bigl[s(a+b)-s_+(c+d)\bigr].
\]
On the other hand,
\[
\begin{aligned}
\mathcal E(T,z)
    &=1-b^2+T^2s(b+c)(b+a),\\
\mathcal E(T,qz)
    &=1-c^2+T^2s_+(c+d)(c+b).
\end{aligned}
\]
Subtracting and factoring gives
\[
\begin{aligned}
\mathcal E(T,qz)-\mathcal E(T,z)
    &=
    b^2-c^2+T^2s_+(c+d)(c+b)\\
    &\qquad -T^2s(b+c)(b+a)\\
    &=
    (b+c)\bigl[
    b-c+T^2s_+(c+d)-T^2s(a+b)
    \bigr]\\
    &=0.
\end{aligned}
\]
The last equality follows directly from the linear identity.
In particular, no division by $b+c$ is required.
Each coefficient of $\mathcal E$ in $T$ is a Laurent
polynomial in $z$. Writing
\[
    \mathcal E(T,z)
    =\sum_{n\geq 0}\sum_{m\in\mathbb Z}
       e_{n,m}(q)T^nz^m,
\]
the identity $\mathcal E(T,qz)=\mathcal E(T,z)$ gives
\[
    (q^m-1)e_{n,m}(q)=0.
\]
Since $q$ is an indeterminate, $q^m-1$ is nonzero in
$\mathbb C(q)$ whenever $m\neq 0$.
Thus $\mathcal E(T,z)$ is independent of $z$.
Finally, if $\eta(T,-1)=1$, then $S(-1)=0$ implies
\[
    \mathcal E(T,z)
    =\mathcal E(T,-1)
    =1-\eta(T,-1)^2
    =0.
\]
\end{proof}

\begin{theorem}\label{thm:higher}
Fix \(q\in\mathbb C^\times\) such that
\(q^m\ne -1\) for every \(m\in\mathbb Z\). The function $\eta(T,z)$ \eqref{GenFuncHigher}, as a formal power series in $T$, is a unique solution of the quadratic TBA equation \eqref{EquationTBA1} with asymptotic $\eta(T,z) = 1 + O(T^2)$.
\end{theorem} 
\begin{proof} Our proof is based on recasting the TBA as a linear equation, together with combinatorial manipulations in the $q$-commutative algebra. We first regard \(q\) as an indeterminate and work over
\(\mathbb C(q)\). At the end of the proof we specialize
to the numerical values of \(q\) allowed in the statement.

\medskip
\noindent
{\em Step 1 - Linearizing the TBA.} Let us use the remarkable fact that the TBA equation can be linearized. It means that there exists a linear difference equation such that any solution of the TBA is necessarily a solution of this linear equation. To see this, we subtract two TBA equations that differ by an elementary shift,
\begin{align}
& \nonumber \big[ 1 - \eta^2(T,z) + T^2 ( \eta( T, q z) + \eta(T, z) )( \eta( T, q^{-1} z) + \eta(T, z) ) ( z + z^{-1} + 2 ) \big] - \textit{}\\
& \big[ 1 - \eta^2(T,qz) + T^2 ( \eta( T, q^2 z) + \eta(T, qz) )( \eta( T, z) + \eta(T, qz) ) ( qz + q^{-1}z^{-1} + 2 ) \big] \mathop{=}^{\star}
\end{align}
and observe that the result is divisible by $\eta(T,z) + \eta(T,qz)$,
\begin{align}
& \nonumber \mathop{=}^{\star} \big[ \eta(T,z) + \eta(T,qz) \big] \times \\
& \left[ \eta(T,qz) - \eta(T,z) + T^2 \frac{(z+1)^2}{z} (\eta( T, q^{-1} z) + \eta(T, z)) - T^2 \frac{(qz+1)^2}{qz} (\eta( T, q^2 z) + \eta(T, qz)) \right].
\end{align}
Factoring out $\eta(T,z) + \eta(T,qz)$, we obtain a linear identity
\begin{align*}
\eta(T,z) - \eta(T,qz) =
T^2\frac{(z+1)^2}{z} (\eta( T, q^{-1} z) + \eta(T, z)) - T^2\frac{(qz+1)^2}{qz} (\eta( T, q^2 z) + \eta(T, qz))
\end{align*}
Introducing $S(z)=2+z+z^{-1}$, one rewrites the equation in the form
\begin{equation}
\eta(T,z) - \eta(T,qz) =
T^2\left[S(z)(\eta( T, q^{-1} z) + \eta(T, z))-S(qz)(\eta( T, q^{2} z) + \eta(T, qz))\right]
    \label{LinearTBAconsequence1}
\end{equation}
% Subtlety here is that not every solution to linear gives a solution to quadratic TBA equation. The structure of quadratic equation suggets the following condition
% \begin{equation*}
%     \eta(T,-1) = \pm 1.
% \end{equation*}

\medskip
\noindent
{\em Step 2 - $\mathcal{A}$-subalgebra.} Now let us introduce 3 elements of $\mathcal{K}$:
\begin{equation}
    A= x+y^{-1},\quad B=x^{-1}+y,\quad Z = xy.
\end{equation}
The quantum plane relation $qxy=yx$ forces the following $A,B$ and $Z$ to form a closed sub-algebra with relations
\begin{align}
    & \nonumber ZA=qAZ, & \,&  AB=2+Z+Z^{-1}=S(Z), \\
    & \nonumber ZB= q^{-1}BZ, & \, &BA=2+qZ+q^{-1}Z^{-1}=S(qZ).
\end{align}
We call this subalgebra $\mathcal{A}$. This subalgebra gives a special case of so-called generalized $q$-Weyl algebra \cite{etingof2021twisted,klyuev2022twisted}. Obviously, the original Laurent polynomial is an element of $\mathcal{A}$-subalgebra
$$
W=x+x^{-1}+y+y^{-1} = A+B \in \mathcal{A},
$$
as well as its powers. Any element of $\mathcal{A}$ reads 
\begin{equation*}\label{eq:arbS}
    U=c(Z)+\sum_i p_i(Z)A^i+\sum_j r_j(Z)B^j,\quad c(Z),p_i(Z),r_j(Z) \in \mathbb{C}(q)[Z,Z^{-1}]
\end{equation*}
The quantum torus carries a \(\mathbb Z\)-grading defined by
\[
\deg x=1,\qquad \deg y=-1,\qquad \deg q=0.
\]
This grading is well defined because the relation \(yx=qxy\)
is homogeneous. In particular,
\[
\deg A=1,\qquad \deg B=-1,\qquad \deg Z=0.
\]
A normally ordered monomial \(x^ay^b\) has degree \(a-b\).
Thus the degree-zero component consists precisely of Laurent
polynomials in \(Z\), since
\[
x^ny^n=q^{-n(n-1)/2}Z^n,\qquad n\in\mathbb Z.
\]
We define \(P_Z\) as the restriction to \(\mathcal A\) of the
linear projection onto the degree-zero component of the
quantum torus. For the normal form displayed above, this gives
\begin{equation}
P_Z[U]=c(Z).
\end{equation}
Indeed, every term \(p_i(Z)A^i\), \(i\geq1\), has positive
degree, whereas every term \(r_j(Z)B^j\), \(j\geq1\), has
negative degree. Hence \(P_Z\) is well defined independently
of the chosen expression for \(U\). Its image is a commutative subalgebra, but need not be central:
the relation \(ZA=qAZ\) must still be respected.
We extend \(P_Z\) coefficientwise to formal power series in \(T\).
Returning back to the initial algebra $\mathcal{K}$ one relates the projector $\langle e_{n,n},\cdot\rangle$ with the projector $\langle Z^n,\cdot\rangle$ as
\begin{equation}
    \langle Z^n, U\rangle = \operatorname{Coeff}_{Z^n} P_Z(U) = q^{{\color{red}{-}}n(n-1)/2}\langle e_{n,n}, U\rangle,\quad \forall U\in \mathcal{A}\subset \mathcal{K},
\end{equation}
due to a simple reordering identity
\begin{equation*}
    Z^n = (xy)^n = q^{-n(n-1)/2}x^ny^n=q^{n(n-1)/2}e_{n,n}.
\end{equation*}
Let us also mention a few important identities for the  projector $P_Z$. For arbitrary $U\in \mathcal{A}$ the following identities hold:
\begin{equation}\label{eq:ABcom}
    P_Z[UA](Z) = P_Z[AU](qZ),\quad P_Z[UB](Z) = P_Z[BU](q^{-1}Z)
\end{equation}
%Indeed let us use the notation of \eqref{eq:arbS}, 
Using the normal form above, we obtain
\begin{equation*}
P_Z[UA] = r_1(Z)S(qZ),\quad P_Z[AU] = P_Z[Ar_1(Z)B]=P_Z[r_1(q^{-1}Z)AB]=r_1(q^{-1}Z)S(Z).
\end{equation*}
A similar proof applies to $P_Z[UB](Z)$.

\medskip
\noindent
{\em Step 3 - Proving TBA.} Introduce the formal inverse 
\begin{equation}\label{eq:R}
    R= (1-TW)^{-1},
\end{equation}
Since $P_Z$ projects to the commutative subalgebra, applying it to $R$ gives just a formal power series in two commuting variables. More importantly, one may relate this series with $\eta$ as follows
\begin{equation}\label{eq:etaR}
    P_Z[R] = \frac{\eta(T,Z)+\eta(T,qZ)}{2}.
\end{equation}
Now let us introduce several elements of the subalgebra generated by $A,B$ and $Z$, whose projections $P_Z$ satisfy nice relations that lead to \eqref{LinearTBAconsequence1} 
\begin{equation}
    H(Z) := P_Z[R],\quad K(Z):=P_Z[AR],\quad J(Z):=P_Z[A^2R],
\end{equation}
which due to the projection all lie in the commutative subalgebra generated by $Z$. Multiplying \eqref{eq:R} by $1-TW$ from the left and from the right, one gets two identities
\begin{equation}\label{eq:Rres}
    R = 1+T(WR),\quad R = 1+T(RW).
\end{equation}
Now multiply both equations by $A$ from the left
\begin{equation*}
    AR = A+T (A^2R + ABR),\quad AR = A + T( ARA + ARB).
\end{equation*}
Applying $P_Z$ to these equations and using the relation \eqref{eq:ABcom}, one gets
\begin{equation*}
    K(Z) = T\left(J(Z) + S(Z)H(Z)\right),\quad K(Z)= T\left(J(qZ) + S(Z)H(q^{-1}Z)\right)
\end{equation*}
These two equations already look very similar to the TBA. Shifting the first equation by $q$ and subtracting the second one, we obtain
\begin{equation*}
    K(qZ)-K(Z) = T\left[S(qZ)H(qZ)-S(Z)H(q^{-1}Z)\right].
\end{equation*}
Combining \eqref{eq:etaR} with the definition of $H$ one gets
\begin{equation*}
    K(qZ)-K(Z) = \frac{T}{2}\left[S(qZ)(\eta( T, q^{2} Z) + \eta(T, 
    qZ))-S(Z)(\eta( T, Z) + \eta(T, q^{-1}Z))\right].
\end{equation*}
In order to show that, $\eta$ solves linear TBA, one needs to prove that 
\begin{equation}\label{eq:Keta}
    K(qZ)-K(Z) = \frac{\eta(T,qZ) - \eta(T,Z)}{2T}.
\end{equation}
To do this, let us introduce involution $\tau$ of $\mathcal{A}$ as follows
\begin{equation*}
    \tau(A)=B,\quad \tau(B)=A,\quad \tau(Z)=Z,\quad \tau(UV)=\tau(V)\tau(U).
\end{equation*}
Since $Z$ is $\tau$-invariant, one gets
\begin{equation*}
    P_Z\circ \tau = P_Z.
\end{equation*}
Moreover, $R$ depends only on $W$, which is also $\tau$-invariant, so one has
\begin{equation*}
    \tau(R)=R.
\end{equation*}
Using these properties one finds another relation for $K(qZ)$ as follows
\begin{equation*}
    P_Z[BR] = P_Z[\tau(BR)]=P_Z[\tau(R)\tau(B)] = P_Z[RA]= K(qZ).
\end{equation*}
To finish the proof, we apply $P_Z$ to \eqref{eq:Rres} which gives
\begin{equation*}
    H(Z) = 1+T (P_Z[AR]+P_Z[BR]) = 1+T(K(Z)+K(qZ)).
\end{equation*}
Using \eqref{eq:etaR} one gets
\begin{equation*}
    \frac{1+\sigma_q}{2} \eta(T,Z) = \frac{1+\sigma_q}{2}(1+2TK(Z)),\quad \sigma_q F(Z)=F(qZ).
\end{equation*}
or equivalently
\begin{equation*}
(1+\sigma_q)\left[\eta(T,Z) - 1 - 2TK(Z)\right]=0.
\end{equation*}
%Assumption of the theorem that $q$ is not a root $-1$, %implies that $1+\sigma_q$ is injective on $\mathbb{C}[q^{\pm %1}][[Z^{\pm 1}]][[T]]$. 
Since \(q^m\neq 1\), action of operator
\[
(1+\sigma_q)Z^m=(1+q^m)Z^m
\]
has a nonzero multiplier for every \(m\in\mathbb Z\).
Thus \(1+\sigma_q\) is injective on
\(\mathbb C(q)[Z,Z^{-1}][[T]]\), which leads to 
\begin{equation*}
    \eta(T,Z) = 1 + 2TK(Z).
\end{equation*}
%This makes \eqref{eq:Keta} an identity. Note at %this point that any solution of the linear TBA %equation with $T^n$ coefficient that is a Laurent %polynomial in $z$ of order $n$ must be unique up %to multiplication by a $z$-independent Taylor %series in $T$ (as a corollary of Section 7, eq. %%(\ref{FundamentalSolutions})). Function %${\widetilde \eta}$ is such, hence up to a $z$-%independent Taylor series in $T$ it is the same as %the unique solution $\eta(T,z)$ of the quadratic %TBA \eqref{EquationTBA1}. Since we already proved %that zero mode of ${\widetilde \eta}$ is the same %as the zero mode of $\eta$, we conclude that %${\widetilde \eta} = \eta$ solves the quadratic %TBA.
Thus the identity \eqref{eq:Keta} holds, and \(\eta\)
satisfies the linear equation \eqref{LinearTBAconsequence1}. We now show that it
also satisfies the quadratic TBA equation. First, we establish the normalization
\[
\eta(T,-1)=1.
\]
Using the normal form above coefficient-wise in \(T\),
write the degree \(-1\) component of the resolvent as
\[
R_{-1}=\rho(T,Z)B,
\qquad
\rho(T,Z)\in\mathbb C(q)[Z,Z^{-1}][[T]].
\]
Since \(A\) has degree \(1\), only this component contributes
to \(P_Z[AR]\). Hence
\[
\begin{aligned}
K(Z)
&=P_Z[AR]
 =A\,\rho(T,Z)B\\
&=\rho(T,q^{-1}Z)\,AB\\
&=\rho(T,q^{-1}Z)S(Z).
\end{aligned}
\]
In particular, \(S(-1)=0\) implies \(K(-1)=0\).
The identity \(\eta(T,Z)=1+2TK(Z)\) therefore gives
\(\eta(T,-1)=1\). 

%Define the residual of the quadratic equation by
%\[
%\begin{aligned}
%\mathcal E(T,z):={}&1-\eta(T,z)^2\\
%&+T^2S(z)
 %\bigl(\eta(T,z)+\eta(T,qz)\bigr)
 %\bigl(\eta(T,z)+\eta(T,q^{-1}z)\bigr).
%\end{aligned}
%\]
%The factorization in (4.4) is an algebraic identity valid
%for arbitrary \(\eta\). Applying it to the present series
%and using the linear equation (4.5), we obtain
%\[
%\mathcal E(T,qz)=\mathcal E(T,z).
%\]
%Each coefficient of \(\mathcal E\) in \(T\) is a Laurent
%polynomial in \(z\) over \(\mathbb C(q)\). If \(e_{n,m}(q)\)
%denotes its coefficient at \(T^nz^m\), the preceding
%identity gives
%\[
%(q^m-1)e_{n,m}(q)=0.
%\]
%Since \(q\) is indeterminate, \(q^m-1\ne0\) for \(m\ne0\).
%Thus \(\mathcal E(T,z)\) is independent of \(z\).
%Evaluating at \(z=-1\), we find
%\[
%\mathcal E(T,z)
%=\mathcal E(T,-1)
%=1-\eta(T,-1)^2
%=0.
%\]
%This proves the quadratic TBA equation over \(\mathbb C(q)\).

The series $\eta$ satisfies the linear identity established
above and the normalization $\eta(T,-1)=1$.
Lemma~\ref{lem:tba-residual} therefore shows that its quadratic
residual $\mathcal E(T,z)$ vanishes.
This proves the quadratic TBA equation over $\mathbb C(q)$.

It remains to specialize \(q\). By the original definition
\eqref{GenFuncHigher}, the coefficient of each power of \(T\) in \(\eta\)
is a finite Laurent polynomial in \(z\), whose coefficients
are rational functions of \(q\). Their possible poles occur
only at \(q=0\) or at zeros of factors \(1+q^m\).
The same is true of the coefficients of \(\mathcal E\).
Consequently, the identity \(\mathcal E=0\) specializes
coefficient-wise to every \(q\in\mathbb C^\times\) satisfying
\(q^m\ne-1\) for all \(m\in\mathbb Z\).

\noindent
Finally, \(P_Z[A]=0\) implies \(K=O(T)\), and therefore
\(\eta=1+O(T^2)\). The uniqueness assertion now follows
from Theorem~\ref{thm:TBA}.
\end{proof}

\section{Method I: from TBA equation to continued fractions}\label{sect:methodI}

In this section, we will derive several compact formulas for the quantum $A$-period in terms of continued fractions. 
We are going to show that equation (\ref{EquationTBA1}) admits a very explicit solution in the class of functions known as \textit{continued fractions}. Given two sequences $a_n,b_n$ for $n=1,2,3,\ldots$, define  the $m$-th \emph{convergent} for $m\geq1$ as
\begin{align}
\mathcal{F}_1[a_n,b_n]=a_1,\quad  \mathcal{F}_m\big[ a_n, b_n \big] = a_1 + \dfrac{b_1}{a_2 + \dfrac{b_2}{\ldots + \dfrac{b_{m-1}}{a_{m}}}},
\end{align}
and the continued fraction as
\begin{equation*}
    \mathcal{F}\big[ a_n, b_n \big] = a_1 + \dfrac{b_1}{a_2 + \dfrac{b_2}{a_3 + \dfrac{b_3}{a_4 + \dfrac{b_4}{a_5 + \ldots }}}}, 
\end{equation*}
Introduce also numerator $\mathcal{N}_m$ and denominator $\mathcal{D}_m$ sequences as follows
\begin{align}
\mathcal{N}_m\big[ a_n, b_n \big] = \left\{ \begin{array}{ll} & 0, \ \ \ m \leq -2 \\ & 1, \ \ \ m = -1 \\ & a_1, \ \ \ m = 0 \\ & a_1 a_2 + b_1, \ \ \ m = 1 \\ & a_{m+1} \mathcal{N}_{m-1}\big[ a_n, b_n \big] + b_{m} \mathcal{N}_{m-2}\big[ a_n, b_n \big], \ \ \ m > 1
\end{array} \right.
\end{align}
\begin{align}
\mathcal{D}_m\big[ a_n, b_n \big] = \left\{ \begin{array}{ll} & 0, \ \ \ m \leq -1 \\ & 1, \ \ \ m = 0 \\ & a_2, \ \ \ m = 1 \\ & a_{m+1} \mathcal{D}_{m-1}\big[ a_n, b_n \big] + b_{m} \mathcal{D}_{m-2}\big[ a_n, b_n \big], \ \ \ m > 1
\end{array} \right.
\end{align}
Then one has
\begin{align}
\mathcal{F}_m\big[ a_n, b_n \big] = \dfrac{\mathcal{N}_{m-1}\big[ a_n, b_n \big]}{\mathcal{D}_{m-1}\big[ a_n, b_n \big]}.
\end{align}
Fix \(q\in\mathbb C^\times\) and set
\[
\mathcal R_z:=\mathbb C[z,z^{-1}],
\qquad
V(z):=z+2+z^{-1}.
\]
We work in the ring \(\mathcal R_z[[T^2]]\) of formal power
series in \(T^2\) whose coefficients are Laurent polynomials
in \(z\). The substitutions \(z\mapsto q^kz\), \(k\in\mathbb Z\),
preserve this ring. Every series with constant term \(1\)
in \(T^2\) is invertible in \(\mathcal R_z[[T^2]]\).

The infinite continued fractions below are understood as
formal limits of their finite convergents: for each \(N\),
the coefficients through \(T^{2N}\) stabilize as the depth
of the convergent increases. The following lemma establishes
the existence of these limits; no analytic convergence
is assumed.
\begin{lemma}[Formal continued fractions]\label{lem:CF}
If $a_n$ is a unit in $\mathcal R_z[[T^2]]$ and $b_n\in T^2\mathcal R_z[[T^2]]$, the
convergents $\mathcal F_m[a_n,b_n]$ have a unique $T$-adic limit.
More precisely, $\mathcal F_{m+1}-\mathcal F_m\in T^{2m}\mathcal R_z[[T^2]]$.
\end{lemma}
\begin{proof}
Every denominator is a unit. Changing the last tail and propagating
outward through $h\mapsto a+b/h$ gains at least one power of $T^2$
at each level, because $b/h-b/\widetilde h
=b(\widetilde h-h)/(h\widetilde h)$.
Comparison at depth $m$ therefore gives a difference in $T^{2m}\mathcal R_z[[T^2]]$.
Completeness gives the limit, independently of any unit terminal tail.
\end{proof}

\begin{theorem}[Continued-fraction solution]\label{thm:CF}
Define the formal fractions
\begin{equation}
f_1(T,z)=\mathcal F[1,-T^2V(q^nz)],\qquad
f_2(T,z)=\mathcal F[1,-T^2V(q^{1-n}z)].
\label{TwoFractions}
\end{equation}
Then the solution of Theorem~\ref{thm:TBA} is
\begin{equation}
\eta(T,z)=\frac{-1-f_1(T,z)+f_2(T,z)}{1-f_1(T,z)-f_2(T,z)}.
\label{TBAContFrac}
\end{equation}
\end{theorem}
\begin{proof}
Lemma~\ref{lem:CF} applies and $f_i=1+O(T^2)$, so the displayed
quotient is defined and equals $1+O(T^2)$. Removing the first level gives
\begin{align}
f_1(z/q)&=1-\frac{T^2V(z)}{f_1(z)},&
 f_2(qz)&=1-\frac{T^2V(qz)}{f_2(z)},\nonumber\\
(f_1(z)-1)f_1(qz)&=-T^2V(qz),&
 (f_2(z)-1)f_2(z/q)&=-T^2V(z).
\label{FractionShifts}
\end{align}
%For the proposed quotient, put $d=f_1+f_2-1$. These identities give
%\[
%\eta(z)+\eta(qz)=\frac2d,\qquad
%T^2V(z)\bigl(\eta(z)+\eta(z/q)\bigr)=\frac{2f_1(1-f_2)}d.
%\]
For the proposed quotient, suppress the dependence on \(T\)
and write
\[
d(z):=f_1(z)+f_2(z)-1,
\qquad
\eta(z)=\frac{2f_1(z)}{d(z)}-1.
\]
The shift identities above imply
\[
f_1(qz)\bigl(f_1(z)-1\bigr)
=
f_2(z)\bigl(f_2(qz)-1\bigr)
=
-T^2V(qz).
\]
Consequently,
\[
\begin{aligned}
f_1(qz)d(z)-f_2(z)d(qz)
&=f_1(qz)\bigl(f_1(z)-1\bigr)
  -f_2(z)\bigl(f_2(qz)-1\bigr)\\
&=0.
\end{aligned}
\]
Since \(d(z)\) and \(d(qz)\) are units, we obtain
\[
\frac{f_1(qz)}{d(qz)}=\frac{f_2(z)}{d(z)}.
\]
It follows that
\[
\begin{aligned}
\eta(z)+\eta(qz)
&=\frac{2f_1(z)}{d(z)}
  +\frac{2f_1(qz)}{d(qz)}-2\\
&=\frac{2\bigl(f_1(z)+f_2(z)\bigr)}{d(z)}-2\\
&=\frac{2}{d(z)}.
\end{aligned}
\]
Replacing \(z\) by \(z/q\) gives
\[
\eta(z)+\eta(z/q)=\frac{2}{d(z/q)},
\qquad
f_1(z)d(z/q)=f_2(z/q)d(z).
\]
Moreover, the defining relation for \(f_2\) gives
\[
T^2V(z)=\bigl(1-f_2(z)\bigr)f_2(z/q).
\]
Combining these identities yields
\[
\begin{aligned}
T^2V(z)\bigl(\eta(z)+\eta(z/q)\bigr)
&=\frac{2T^2V(z)}{d(z/q)}\\
&=\frac{2f_1(z)T^2V(z)}{f_2(z/q)d(z)}\\
&=\frac{2f_1(z)\bigl(1-f_2(z)\bigr)}{d(z)}.
\end{aligned}
\]
%Their product equals $\eta(z)^2-1=4f_1(1-f_2)/d^2$, proving
%\eqref{EquationTBA}. All these identities can be checked after
%clearing unit denominators; no inversion of $f_i-1$ is required.
%Uniqueness finishes the proof.
Multiplying the two identities just proved, we obtain
\[
\begin{aligned}
&T^2V(z)
 \bigl(\eta(z)+\eta(qz)\bigr)
 \bigl(\eta(z)+\eta(z/q)\bigr)\\
&\qquad=\frac{4f_1(z)\bigl(1-f_2(z)\bigr)}{d(z)^2}
 =\eta(z)^2-1,
\end{aligned}
\]
where the last equality follows from
\[
\bigl(1+f_1(z)-f_2(z)\bigr)^2-d(z)^2
=4f_1(z)\bigl(1-f_2(z)\bigr).
\]
Thus the proposed expression satisfies the quadratic
TBA equation~\eqref{EquationTBA1}.

All divisions used above are by formal series with constant
term \(1\); in particular, no division by \(f_i-1\) is needed.
Moreover, \(f_1,f_2=1+O(T^2)\) implies
\(\eta=1+O(T^2)\).
By the uniqueness statement of Theorem~\ref{thm:TBA}, this is the
required formal solution.
\end{proof}

\section{The barrier method}\label{sec:barrier}

Here is a second proof of Theorem~\ref{thm:higher}.
Instead of the linearization \eq{LinearTBAconsequence1}, we expand $(1-TW)^{-1}$ into walks and cut each walk where it returns to the diagonal; this gives the formula \eq{TBAContFrac} for $\eta$ directly.
Its constant term gives Theorem~\ref{thm:zero}.

\paragraph{The picture.} Let us describe the intuition behind the method first, and give a formal proof afterwards.
Expanding $(1-TW)^{-1}=\sum_nT^nW^n$, with $W=x+x^{-1}+y+y^{-1}$, gives all the words in the letters $x^{\pm1}$, $y^{\pm1}$.
We read a word from right to left, as an operator acting on a vector, and draw it as a walk on the grid $\mathbb{Z}^2$: its rightmost letter is the first step, and the letters $x^{\pm1}$ and $y^{\pm1}$ are the steps $\pm(1,0)$ and $\pm(0,1)$.
By $yx=qxy$, the word of a walk from $0$ to $(a,b)$ is $q^{c+ab/2}x^ay^b$, where $c$ is the area enclosed by the walk and the segment from its end back to $0$, counted positively when this loop runs counterclockwise; so the series below count walks with their areas.
We use the notation of Step~2 of the proof of Theorem~\ref{thm:higher}: the grading $\deg x=1$, $\deg y=-1$ of $\mathcal{K}$, the elements $A=x+y^{-1}$ and $B=x^{-1}+y$ of degrees $1$ and $-1$, so that $W=A+B$, and $Z=xy$.
The \emph{height} of a point $(a,b)$ is $a-b$, the degree of $x^ay^b$: a step of $A$ ($x$ or $y^{-1}$, right or down) raises it by one, and a step of $B$ ($x^{-1}$ or $y$, left or up) lowers it by one.
The method takes its name from the following intuitive picture: the diagonal through $0$ is a barrier that registers every visit of the walk.
Every step changes the height by $\pm1$, so, by continuity, a walk cannot cross the barrier without landing on it.

\begin{figure}[ht]
\centering
\begin{tikzpicture}[scale=0.7]
  \draw[step=1,gray!30,very thin] (-1.5,-0.5) grid (5.5,5.5);
  \draw[dashed,gray] (-0.5,-0.5) -- (5.5,5.5);
  \draw[very thick,red,->] (0,0)--(0,1)--(1,1);
  \draw[very thick,blue,->] (1,1)--(2,1)--(3,1)--(3,2)--(3,3);
  \draw[very thick,red,->] (3,3)--(2,3)--(2,4)--(3,4)--(4,4);
  \foreach \p in {(0,0),(1,1),(3,3),(4,4)} \fill \p circle (3pt);
  \node[red] at (-0.5,1.2) {$-$};
  \node[blue] at (3.7,1.4) {$+$};
  \node[red] at (1.6,4.2) {$-$};
  \node[below right] at (0,0) {$0$};
\end{tikzpicture}
\qquad
\begin{tikzpicture}[scale=0.8]
  \fill[gray!20] (1,0) -- (3,2) -- (3,3) -- (1,1) -- cycle;
  \draw[step=1,gray!30,very thin] (-0.5,-1.5) grid (4.5,4.5);
  \draw[dashed,gray] (-0.5,-0.5) -- (4.5,4.5);
  \draw[dotted,thick,gray] (-0.5,-1.5) -- (4.5,3.5);
  \draw[very thick,->] (0,0)--(1,0);
  \draw[very thick,->] (0,0)--(0,-1);
  \draw[thick,decorate,decoration={snake,amplitude=1pt,segment length=6pt}] (1,0)--(3,2);
  \draw[very thick,->] (3,2)--(3,3);
  \draw[very thick,->] (3,2)--(2,2);
  \fill (0,0) circle (3pt);
  \node at (2.45,0.6) {$w$};
  \node at (0.35,-0.35) {$A$};
  \node at (2.7,2.35) {$B$};
\end{tikzpicture}
\caption{Left: the diagonal (dashed) is the barrier; a walk from $0$ to the diagonal is cut at its visits to the barrier into primitive walks, positive ($+$) below the barrier and negative ($-$) above it.
Right: a primitive positive walk is a step right or down (a letter of $A=x+y^{-1}$) onto the dotted line one level below the barrier, a positive walk $w$ on or below the dotted line, and a step left or up (a letter of $B=x^{-1}+y$) back to the barrier; its word is $B\,w\,A$.
Compared with running along the barrier, running along the dotted line adds the shaded strip to the area: $k$ if $w$ runs from $(a,b)$ to $(a+k,b+k)$, here $k=2$.}
\label{fig:bar-barrier}
\end{figure}
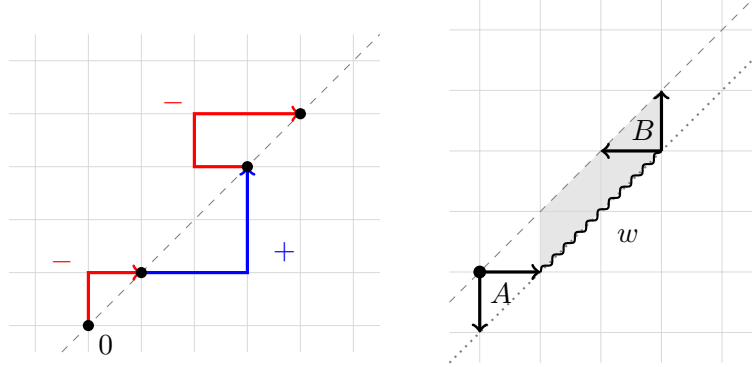

The elements of degree zero are the Laurent polynomials in $Z$.
We write $P_Z$ for the projection of $\mathcal{K}$ onto its component of degree zero (Step~2 restricts it to $\mathcal{A}$), and extend it coefficientwise to series in $T$.
A word has degree zero if and only if its walk ends at the height of its start, that is, on the diagonal through its start.
The projection $P_Z$ keeps exactly these words, so $H=P_Z[(1-TW)^{-1}]$ is the generating function of the walks from $0$ that end on the diagonal $\{(m,m)\}$, in which every walk with $n$ steps is counted with $T^n$ times its word.
That word is $q^{c+m/2}Z^m$, with $c$ the area of the walk, so $H(q^{-1/2}z)$ is the sum of $T^nq^cz^m$ over these walks.
Such a walk is \emph{primitive} if it is non-empty and returns to the height of its start only at its end, and \emph{positive} (resp.\ \emph{negative}) if it never goes below (resp.\ above) the height of its start.
The visits to the barrier cut every walk from $0$ to the diagonal, uniquely, into primitive walks, each of them positive (below the barrier) or negative (above it); see Figure~\ref{fig:bar-barrier}, left.
Let $H^p$, $H_\pm$ and $H^p_\pm$ be the generating functions of these subclasses of walks, where the superscript $^p$ selects only the primitive walks and the subscript $_\pm$ selects only the positive or negative walks; like $H$, they have degree zero.

Consider a primitive positive walk.
Its first step has to raise the height (come from $A$), its last step has to lower it (come from $B$), and the steps from the second to the penultimate form a positive, but not necessarily primitive, walk that has two fewer steps; conversely, these three pieces can be chosen independently (Figure~\ref{fig:bar-barrier}, right).
The middle piece runs on or below the next line, of height one, which serves as its barrier.
So the positive walks are described by the same walks one level lower, and so on; this gives a continued fraction, with one level of the fraction for each line.
Running along the next line instead of along the barrier adds to the area the strip between the two lines, whose area is $k$ for a run from $(a,b)$ to $(a+k,b+k)$; this is why the level $h$ of the continued fraction involves $q^h$.

\paragraph{The proof.}
By Step~2 of the proof of Theorem~\ref{thm:higher}, $BA=S(qZ)$ and $AB=S(Z)$, where $S(Z)=Z+2+Z^{-1}$ is the $V(Z)$ of \eq{TwoFractions}.
For any $Q$ of degree $k$ and any $u(Z)$ of degree zero one has
\beq\label{eq:bar-swap}
u(Z)\,Q=Q\,u(q^kZ),\qquad Q\,u(Z)=u(q^{-k}Z)\,Q.
\eeq
Indeed, the equation is bilinear in $u$ and $Q$, so it suffices to check it on monomials, where it follows from $Zx=qxZ$ and $Zy=q^{-1}yZ$, that is, from $yx=qxy$.
By \eq{eq:bar-swap}, for any $u(Z)$ of degree zero one has
\beq\label{eq:bar-AHB}
B\,u(Z)\,A=S(qZ)\,u(qZ),\qquad A\,u(Z)\,B=S(Z)\,u(q^{-1}Z).
\eeq
For $u=Z^k$ the first identity reads $B\,Z^kA=q^kZ^k\,BA$; its exponent $k$ is the area of the strip in the picture (shaded in Figure~\ref{fig:bar-barrier} for $k=2$).

Algebraically, the decomposition of a primitive positive walk described in the picture is the following relation between generating functions in the non-commutative algebra $\mathcal{K}$:
\beq\label{eq:bar-peel}
H^p_+=T^2\,B\,H_+\,A,
\eeq
which thanks to \eq{eq:bar-AHB} can be transformed into an equation in the commutative algebra of series in $T$ whose coefficients are Laurent polynomials in $q$ and $Z$:
\beq\label{eq:bar-Hp-phi}
H^p_+(Z)=T^2S(qZ)\,H_+(qZ).
\eeq
Any positive walk can be uniquely decomposed into a sequence of primitive positive walks, so
\beq\label{eq:bar-first-return-plus}
H_+=1+H^p_++(H^p_+)^2+\dots=\frac{1}{1-H^p_+},
\eeq
and in \eq{eq:bar-Hp-phi} one can substitute $H_+(qZ)$ by $1/(1-H^p_+(qZ))$ to obtain the following recursion, and hence a continued fraction, for $H^p_+$:
\beq\label{eq:bar-contfrac}
H^p_+(Z)=\frac{T^2S(qZ)}{1-H^p_+(qZ)}=\cfrac{T^2S(qZ)}{1-\cfrac{T^2S(q^2Z)}{1-\cfrac{T^2S(q^3Z)}{1-\dotsb}}}\,.
\eeq
Here the continued fraction is the $T$-adic limit of Lemma~\ref{lem:CF}, which does not depend on the unit tails $1-H^p_+(q^kZ)$; comparing it with the definition \eq{TwoFractions} of $f_1$, we find $1-H^p_+=f_1$ and $H_+=1/f_1$.
Before \eq{eq:bar-AHB} is used, \eq{eq:bar-peel} and \eq{eq:bar-first-return-plus} already give $H_+$ as a continued fraction in $\mathcal{K}$: it is the image of Flajolet's continued fraction for paths in non-commuting variables \cite[Theorem~1]{Flajolet1980} when every rise is replaced by $TA$, every fall by $TB$ and every level step by $0$, and words are read from right to left.

The automorphism $x\mapsto x^{-1}$, $y\mapsto y^{-1}$ of $\mathcal{K}$ rotates the grid by $\pi$.
It fixes $W$, exchanges $A$ and $B$, sends $Z$ to $x^{-1}y^{-1}=q^{-1}Z^{-1}$, and multiplies all heights by $-1$, thus swapping positive and negative walks, and positive primitive walks with negative primitive walks.
This implies that
\beq\label{eq:bar-rotation}
H^p_-(Z)=H^p_+(q^{-1}Z^{-1}),
\eeq
so that $1-H^p_-=f_2$: since $S(Z^{-1})=S(Z)$, the substitution $Z\mapsto q^{-1}Z^{-1}$ interchanges $f_1$ and $f_2$.

Since every primitive walk is positive or negative, the generating function of all primitive walks equals
\beq\label{eq:bar-Hp-sum}
H^p=H^p_++H^p_-=2-f_1-f_2.
\eeq
Finally, the cut at the barrier gives
\beq\label{eq:bar-H}
H=\frac{1}{1-H^p}=\frac{1}{f_1+f_2-1}.
\eeq

\paragraph{The TBA.}
By \eq{eq:etaR}, $H=\frac12\bigl(\eta(T,Z)+\eta(T,qZ)\bigr)$.
The proof of Theorem~\ref{thm:CF} shows that the right side of \eq{eq:bar-H} is $\frac12\bigl(\tilde\eta(Z)+\tilde\eta(qZ)\bigr)$, where $\tilde\eta$ is the right side of \eq{TBAContFrac}.
Hence $(1+\sigma_q)(\eta-\tilde\eta)=0$, where $\sigma_qu(Z)=u(qZ)$; under the hypothesis $q^m\neq-1$ of Theorem~\ref{thm:higher}, $1+\sigma_q$ is injective (it multiplies $Z^m$ by $1+q^m\neq0$, as in Step~3 of its proof), so $\eta=\tilde\eta$.
So the function $\eta$ of \eq{GenFuncHigher} is given by \eq{TBAContFrac}, and by Theorem~\ref{thm:CF} it is the unique solution of the TBA equation \eq{EquationTBA1} given by Theorem~\ref{thm:TBA}; this proves Theorem~\ref{thm:higher} again, without the linearization.

\paragraph{The zero mode.}
Taking constant terms in \eq{eq:bar-H} gives the quantum period.

\begin{theorem}[Zero mode and the quantum period]\label{thm:zero}
%Without any use of Theorem~\ref{thm:higher},
\begin{equation}
G(T)=\CT_z\eta(T,z)
=\CT_z\frac{-1-f_1(T,z)+f_2(T,z)}{1-f_1(T,z)-f_2(T,z)}.
\label{PeriodContFrac1}
\end{equation}
Equivalently, $\CT_z\eta=\theta_\zeta\Pi_A(\zeta;q)$ at $\zeta=T^2$,
where $\theta_\zeta=\zeta\,\partial/\partial\zeta$.
\end{theorem}
\begin{proof}
Here $\eta$ denotes the solution given by Theorem~\ref{thm:CF}, the $\tilde\eta$ above, and $q\in\mathbb{C}^\times$ is arbitrary.
Since $\CT_ZP_Z[L]=\langle e_{0,0},L\rangle$ for every $L\in\mathcal{K}$, \eq{eq:bar-H} gives $G(T)=\CT_ZH=\CT_z\,d(z)^{-1}$, where $d=f_1+f_2-1$.
The proof of Theorem~\ref{thm:CF} shows that $\eta(z)+\eta(qz)=2/d(z)$; since the substitution $z\mapsto qz$ preserves constant terms, $\CT_z\eta=\CT_z\,d^{-1}=G(T)$.
Neither step uses Theorem~\ref{thm:higher}.

For the final statement, use the normalization
\[
\Pi_A(\zeta;q)
=\log\zeta-2F(\sqrt{\zeta})+C(q),
\qquad \zeta=T^2,
\]
where \(C(q)\) is independent of \(\zeta\). Then
\[
\theta_\zeta\Pi_A
=1-TF'(T)
=G(T)
=\operatorname{CT}_z\eta(T,z).
\]
\end{proof}

\noindent
This is the relation (2.10) of \cite{HuangTBA} for the zero mode, in our normalization.
There it is stated for a general mass parameter and checked by expansion in the modulus, without proof; above, it is proved for our $W$ as an identity of series in $T$.
%, without Theorem~\ref{thm:higher}.
%\input{barrier_section}

\section{Method II: from linearized TBA to fast algorithms}\label{sect:methodII}

In this section we derive another representation for the A-period -- in terms of an explicit recursion that provides an efficient way to compute it order by order. We will use the remarkable fact that was already mentioned and proven --- TBA equation can be linearized to \eqref{LinearTBAconsequence1}, which reads as
\begin{align}
\eta(T,z) - \eta(T,qz) =
T^2\frac{(z+1)^2}{z} (\eta( T, q^{-1} z) + \eta(T, z)) - T^2\frac{(qz+1)^2}{qz} (\eta( T, q^2 z) + \eta(T, qz))
\label{LinearTBAconsequence}
\end{align}
Define now an operator
\begin{align}
{\widehat A} f(z) := \int\limits_{-1}^{z} \dfrac{d_q y}{y} \left[ \frac{(y+1)^2}{y} (f( q^{-1} y) + f(y)) - \frac{(qy+1)^2}{qy} (f(q^2 y) + f( qy)) \right]
\end{align}
where the $q$-integral is defined on Laurent polynomials in $z$ without constant term
\begin{align}
 \int\limits_{-1}^{z} \dfrac{d_q y}{y} y^n = \dfrac{z^{n}-(-1)^n}{1-q^n}, \ \ \ \forall n \neq 0
\end{align}
Note that the action of ${\widehat A}$ can be made completely explicit
\begin{align}
{\widehat A} \ \sum\limits_{i} c_{|i|} z^i := \sum\limits_{i} c^{\prime}_{|i|} z^i
\label{Algorithm1}
\end{align}
where
\begin{align}
c^{\prime}_{0} = 4 c_0 + 2 \dfrac{q^2+3q+1}{q} c_1 + \sum\limits_{i = 2}^{\infty} \left( \dfrac{2(-1)^i(1+q^i)^2(1-q)^2}{q(1-q^{i-1})(1-q^{i+1})} \right) c_i,
\end{align}
\begin{align}
c^{\prime}_{j} = \dfrac{(1-q^{2j-1})(1+q^{j-1})}{q^{j-1}(1-q^{j})} c_{j-1} + 2 \dfrac{(1+q^j)^2}{q^j} c_{j} +  \dfrac{(1-q^{2j+1})(1+q^{j+1})}{q^{j+1}(1-q^{j})} c_{j+1}, \ \ \ j \geq 1
\end{align}
Then, the linearized TBA equation takes the form
\begin{align}
-1 + \eta(T,z) =
T^2 {\widehat A} \eta(T,z)
\end{align}
and its solution is
\begin{align}
\eta(T,z) = \dfrac{1}{1 - T^2 {\widehat A}} \cdot 1
\label{Algorithm2}
\end{align}
This gives an explicit algorithm to compute the coefficients of $\eta(T,z)$ up to $T^{2n}$ in $O(n^2)$ operations. Indeed, every next order in $T^2$ is obtained from the previous one by the action of the same operator ${\widehat A}$, and it has $O(n)$ coefficients (that each take a fixed number of operations to compute). As to compute the $n$-th order we first need to compute all the previous ones, we see this can be done in $O(n^2)$ time.

What is even more remarkable is that the following algorithm computes $1/G(T)$. Let ${\widehat B}$ be the following operator on symmetric Laurent polynomials in $z$
\begin{align}
{\widehat B} \ \sum\limits_{i} c_{|i|} z^i := \sum\limits_{i} c^{\prime}_{|i|} z^i
\end{align}
where
\begin{align}
c^{\prime}_{0} = 0,
\end{align}
\begin{align}
c^{\prime}_{j} = \dfrac{(1-q^{2j-1})(1+q^{j-1})}{q^{j-1}(1-q^{j})} c_{j-1} + 2 \dfrac{(1+q^j)^2}{q^j} c_{j} +  \dfrac{(1-q^{2j+1})(1+q^{j+1})}{q^{j+1}(1-q^{j})} c_{j+1}, \ \ \ j \geq 1
\end{align}
Then,
\begin{align}
\dfrac{\eta(T,z)}{G(T)} = \dfrac{1}{1 - T^2 {\widehat B}} \cdot 1
\end{align}
This follows at once from the fact that $\dfrac{1}{1 - T^2 {\widehat B}} \cdot 1$ and $\dfrac{1}{1 - T^2 {\widehat A}} \cdot 1$ differ by normalization with some overall function of $T$, and the constant term in $z$ of both the l.h.s. and the r.h.s. is 1. Now we can specialize $z = -1$ and, since from the TBA
\begin{align}
\eta(T,-1) = 1
\end{align}
we finally obtain
\begin{align}
1/G(T) = \dfrac{1}{1 - T^2 {\widehat B}} \cdot 1 \ \ \ \Big|_{z=-1}
\end{align}
Expanding the r.h.s. here, we find an explicit combinatorial formula for $1/G(T)$
\begin{align}
\dfrac{1}{G(T)} = 1 + 2 \sum\limits_{n = 1}^{\infty} \ \ T^{2n} \ \ \left(  \sum\limits_{k_1, \ldots, k_n \in \{-1,0,1\}} \ \prod\limits_{i = 1}^{n} C^{(k_i)}_{k_1 + \ldots + k_{i-1}}
\right)
\end{align}
where
\begin{align}
C^{(-1)}_{j} &=
\begin{cases}
0, & j \leq 1, \\[4pt]
-\dfrac{(1-q^{2j-1})(1+q^{j})}{q^{j}(1-q^{j-1})}, & j>1,
\end{cases}
\\
C^{(0)}_{j} &=
\begin{cases}
0, & j \leq 0, \\[4pt]
2\dfrac{(1+q^j)^2}{q^j}, & j>0,
\end{cases}
\\
C^{(+1)}_{j} &=
\begin{cases}
0, & j<0, \\[4pt]
-\dfrac{(1-q^{2j+1})(1+q^{j})}{q^{j}(1-q^{j+1})}, & j\geq0.
\end{cases}
\end{align}
This result can be compared to the explicit combinatorial formula of \cite{Gan:2024rzb} for the coefficients of $T^n$ in $G(T)$ itself, which appears to be much more complicated. It is quite surprising that our result for $1/G(T)$ has this much more compact form, and it would be interesting to understand fully the mechanisms behind this observation, as well as its implications.

\section{Method III: from linearized TBA to orthogonal polynomials}\label{sect:methodIII}
Keep the generic-$q$ convention of Method II. Define
$C_k(T)=[z^k]\eta(T,z)=[z^{-k}]\eta(T,z)$ for $k\geq0$, so that
\begin{equation}
\eta(T,z)=C_0(T)+\sum_{k\geq1}(z^k+z^{-k})C_k(T),\qquad C_0(T)=G(T).
\end{equation}
Set $u=T^2$ and $\lambda=u^{-1}$. We regard $C_k$ as a series in
$u$; $K((u))$ denotes formal Laurent series with only finitely many
negative powers. The scalar recurrence coefficients $A_k,B_k$ below
are distinct from the operators $\widehat A,\widehat B$.
The nonzero-mode recurrence becomes
\begin{equation}
C_k=A_{k-1}C_{k-1}+B_{k-2}C_{k-2},\qquad k\geq2,
\label{ThreeTerm}
\end{equation}
where, for integer $k\geq-1$ (the index $-1$ is used for a tail below),
\begin{align}
A_k&=-q\frac{\bigl(2(1+q^k)^2-\lambda q^k\bigr)(1-q^k)}
                  {(1+q^{k+1})(1-q^{2k+1})},\label{Ak}\\
B_k&=-q^2\frac{(1-q^{2k+1})(1+q^k)}{(1+q^{k+2})(1-q^{2k+3})}.
\label{Bk}
\end{align}
In particular $A_0=0$ and $B_{-1}=1$. The polynomial variable is
$\lambda=T^{-2}$.
For $k\geq1$, $A_k$ has degree one in $\lambda$ with nonzero leading
coefficient, and $B_k\ne0$ in $K$.

Put $N_j=\mathcal N_j[A_n,B_n]$ and $D_j=\mathcal D_j[A_n,B_n]$,
with the same initial conditions as in Method I, and define
\begin{equation}
\omega(u)=\mathcal F[A_{n-1},B_{n-1}]
=\frac{B_0}{A_1+\displaystyle\frac{B_1}{A_2+\cdots}}.
\label{OmegaFraction}
\end{equation}
Although $A_k$ contains $u^{-1}$, this fraction is well defined:
if $r_k=(A_k+B_kr_{k+1})^{-1}$ for $k\geq1$, then
\[
r_k=u\bigl(uA_k+uB_kr_{k+1}\bigr)^{-1}\in uK[[u]].
\]
The denominator has nonzero constant term. A change of tail gains two
powers of $u$ at each step, proving stabilization; $\omega=B_0r_1$.
The solution selected by the support bound $C_k=O(u^k)$ satisfies
$C_1/C_0=-\omega$, by backward substitution in \eqref{ThreeTerm}.
Indeed its leading coefficients are nonzero for generic $q$ (specialize
the TBA recursion to $q=1$ to see this), so the ratios used in this
backward substitution exist in $K((u))$ and are $O(u)$.
The fundamental-solution formula is consequently
\begin{equation}
\frac{C_k}{C_0}=B_0D_{k-2}-\omega N_{k-2},\qquad k\geq1,
\label{FundamentalSolutions}
\end{equation}
where $D_{-1}=0$, $N_{-1}=1$. Its first two cases are
$-\omega$ and $B_0-A_1\omega$, and induction proves the rest.

\paragraph{Formal orthogonality.}
For a formal Laurent differential we use
\[
\Res_{u=0}g(u)\,du=[u^{-1}]g(u),\qquad
\Res_{u=0}g(u)\frac{du}{u}=[u^0]g(u),
\]
with no factor $2\pi i$.
Define the $K$-linear functionals on $K[\lambda]$, for $i=1,2$, by
\begin{align}
\mathcal L_i[p]&=\Res_{u=0}\Omega_i(u)p(u^{-1})\frac{du}{u},\\
\Omega_1(u)&=\frac{\omega(u)}u,\qquad
\Omega_2(u)=\frac1{u\omega(u)}+\frac1{2u^2}-\frac{(q+1)^2}{qu}.
\label{ResidueWeights}
\end{align}
The negative powers in $\Omega_2$ cancel. Equivalently,
$\Omega_2=u^{-1}\mathcal F[A_{n-2},B_{n-2}]$; the first two levels
reduce it to the expression above. These are formal residue
functionals, without any assertion of positivity or contour convergence.
They obey
\begin{align}
\mathcal L_1[N_iN_j]&=
-\frac{2q^{i+1}(1-q^2)}{1-q^{2i+4}}\delta_{ij},&&i,j\geq-1,
\label{NormN}\\
\mathcal L_2[D_iD_j]&=
\frac{q^{i-3}(1-q^4)(1+q^2)(1+q+q^2)^2}{2(1-q^{2i+4})}\delta_{ij},
&&i,j\geq0.
\label{NormD}
\end{align}
For completeness, after making the polynomials monic, the recurrences
are of the form $p_{r+1}=(\lambda-\alpha_r)p_r-\beta_rp_{r-1}$
with $\beta_r\ne0$ for $r\geq1$. Here monic means leading coefficient
one, and $\alpha_r,\beta_r$ are the scalar coefficients obtained by
rescaling the displayed three-term recurrences. The corresponding
Jacobi fraction is the expansion at $\lambda=\infty$ of
\[
\mathcal J(\lambda)=
\frac{1}{\lambda-\alpha_0-\displaystyle\frac{\beta_1}
{\lambda-\alpha_1-\displaystyle\frac{\beta_2}{\lambda-\alpha_2-\cdots}}}
=\sum_{n\geq0}\mu_n\lambda^{-n-1}.
\]
It defines the normalized moment functional by
$\mathcal L[\lambda^n]=\mu_n$; a different value of $\mathcal L[1]$
rescales all moments. In this construction, multiplying its
convergents by their denominators proves orthogonality by comparison
of Laurent coefficients. Applying this construction to the two tails
in \eqref{ResidueWeights} gives the two functionals above. More
explicitly, let $a_k=[\lambda]A_k$ denote the coefficient of $\lambda$
in $A_k$. In the $N$ family write $h_i=\mathcal L_1[N_i^2]$, and in
the $D$ family write $h_i=\mathcal L_2[D_i^2]$. These ``squared norms''
refer to the bilinear pairing $(p,r)\mapsto\mathcal L_j[pr]$, without
complex conjugation or positivity. They satisfy
$h_i/h_{i-1}=-B_{i+1}a_{i+1}/a_{i+2}$.
This holds for $i\geq0$ in the $N$ family and $i\geq1$ in the $D$ family.
Their initial values are $h_{-1}=-2$ for $N$ and
$h_0=q^{-3}(1+q^2)(1+q+q^2)^2/2$ for $D$; substitution proves
\eqref{NormN}--\eqref{NormD}. The degrees of $N_i,D_i$ are $i+1,i$,
respectively. A functional $\mathcal L$ is called \emph{quasi-definite}
if every moment matrix $(\mathcal L[\lambda^{a+b}])_{0\leq a,b\leq n}$
has nonzero determinant. The complete polynomial families and nonzero
norms above prove this property for generic $q$. Finally, $\eta(T,-1)=1$ gives
\begin{align}
1=C_0(T)+2\sum_{k\geq1}(-1)^kC_k(T),
\label{OrthogonalPeriod1}
\end{align}
\begin{align}
\frac1{G(T)}=1+2\sum_{k\geq1}(-1)^k
\bigl(B_0D_{k-2}-\omega N_{k-2}\bigr). 
\label{OrthogonalPeriod2}
\end{align}
The sum is formal: each bracket is $C_k/C_0=O(u^k)$. The two terms
inside a bracket must be combined before summing, since they contain
cancelling negative powers. This yields the quantum period from the
same three-term recurrence.

\section{Picard-Fuchs, perturbative differential and difference equations}\label{sect:PF}

In this section we discuss in more detail how one could try to extend the classical differential equation for the periods to the quantum level.  
For this it is convenient to consider the function $G(T)$ defined in \eq{QuantumAPeriodG}. As it contains only even powers of $T$, we denote
\beq
    \tau=T^2 \ .
\eeq
One can directly show that for $q=1$ the function $G$  satisfies the following differential equation
\beq
\label{ode1}
   \left(16 \tau^2-\tau\right) G''(\tau)+(32 \tau-1) G'(\tau)+4 G(\tau)=0
\eeq
% \beq
%     \left.(\theta_T^2-16T^2(\theta_T+1)^2)G(T)\right|_{q=1}=0
% \eeq
which can be viewed as the Picard-Fuchs equation for the classical period. This is a manifestation of a more general phenomenon discussed in the previous section.

A natural question is how this equation is deformed when we introduce the noncommutativity parameter $q$. Surprisingly, it seems that there is no simple answer to this question. What we observed explicitly is that we can find differential equations for $G(\tau)$ order by order in the expansion around $q=1$, setting
\beq
    q=e^\epsilon
\eeq
with $\epsilon\to 0 $. 
The order of the differential equation (the number of derivatives) grows with the order of this expansion. In particular, the leading correction reads
\beqa && \label{ode3}
    \frac{1}{768} \tau \epsilon ^2 G^{(3)}+\left(\tau \left(-\frac{7 \epsilon ^2}{192}-\frac{1}{4}\right)+4 \tau^2+\frac{\epsilon
   ^2}{384}\right) G'' \\ \nn && +\left(\frac{4 \tau^2 \epsilon ^2}{3}+\tau \left(8-\frac{\epsilon ^2}{12}\right)-\frac{\epsilon
   ^2}{24}-\frac{1}{4}\right) G'+G\left(\frac{2 \tau \epsilon ^2}{3}-\frac{\epsilon ^2}{48}+1\right)=0
\eeqa
which holds up to $O(\epsilon^3)$ terms.
Note there are ambiguities in writing this equation (one can always e.g. add the 0th order equation or its derivatives), we tried to fix them to keep the equation as simple as possible. Notice there is no $O(\epsilon)$ correction as the function $G$ is symmetric under $q\to 1/q$. We have computed several more orders of the expansion but did not find any simple structure\footnote{Alternatively, one can also try to fix the order of the equation to be 2, while having it reduce at $\epsilon\to 0$ not literally to the original P-F equation, but to the P-F equation multiplied by some function. We find that this is possible (with polynomial coefficients in the equation), but at $\epsilon\to 0$ the equation then reduces to P-F equation times $(16\tau-1)^n$ where $n$ grows with the imposed order of the expansion in $\epsilon$. In other words, the equation does not have a `stable' limit to all orders in $\epsilon$.}.  We give the equation to higher order in appendix \ref{app:ode}.  We note that related observations were explored in \cite{HuangTBA,Huang2012NS} (see also \cite{HuangEtAl2015,HuangSugimotoWang2020,MironovMorozov2010SU2,AganagicEtAl2012}), where in particular higher orders in $\epsilon$ in $G$ were shown to be given by certain differential operators acting on the 0th order. 

Naively, one may expect that these perturbative corrections could resum into simple finite-difference operators, providing a "quantisation" of the original P-F equation. However, this does not seem to be the case. We attempted to find a finite-difference equation of the form
\beq  \label{dieq} \sum_{k=0}^{K} P_k(\tau)G(q^k\tau)=0
\eeq
where $P_k(\tau)$ are polynomials in $\tau$. By plugging in the series expansion of $G$ in $\tau$ (at finite $q$) and imposing \eq{dieq}, we found that no such difference equation exists of order $K\leq 8$ with polynomial coefficients of degree up to $18$.

It would be remarkable if a proper resummation of all orders in the expansion around $q=1$ could somehow be done. Furthermore this would provide a highly efficient way to compute $G(\tau)$ as a series in $\tau$. In the absence of such a resummation, in the next sections we instead explore the TBA equations for our function and find a number of new structures there.

% \section{Conclusions and open problems}

%polynomial $W = x + x^{-1} + y + y^{-1}$ in it. The canonical trace identifies the %quantum period with an area-weighted
%admits a continued-fraction expression, whose zero mode equals the period by %Theorem~\ref{thm:zero}.
%Its linearization supplies a quadratic-cost coefficient algorithm and
%formal orthogonality functionals. The full higher-mode realization is proved, while %the single-fraction logarithmic identity remains a conjecture.
%
%For practical computation of the quantum period, especially important is the %algorithm eqs. (\ref{Algorithm1})-(\ref{Algorithm2}). In computations, we see that %this algorithm is significantly faster, than all other methods. It computes the %quantum period up to $T^{2n}$ in $O(n^2)$ operations. This can be compared, e.g. to %the naive algorithm that follows from expanding the power of noncommutative Laurent %polynomial in eq. (\ref{QuantumAPeriod}) -- there, we would get $O(n^3)$ %operations. 
%
%We leave extensions to other mirror polynomials, efficient computation of the zero %mode alone, and comparison with the canonical spectral polynomials
%as natural directions for further work.

\section{Conclusions and open problems}\label{sect:conclusions}

We have established an algebraic correspondence between all
Laurent modes of the formal TBA solution and the normalized
diagonal coefficients of the quantum-torus resolvent
$(1-TW)^{-1}$, for
$W=x+x^{-1}+y+y^{-1}$ with $yx=qxy$.
In particular,
\[
\operatorname{CT}_z\eta(T,z)
=
G(T)
=
1-T\frac{dF}{dT},
\qquad F(0)=0,
\]
so the quantum A-period is recovered from the zero mode
by formal integration.
Continued fractions, lattice-path decompositions and
the linearized functional equation provide complementary
descriptions of this correspondence, leading to explicit
coefficient recursions and formal orthogonal polynomials.
Several questions remain open.

\paragraph{Scalar $q$-difference equations in the modulus.}
The linearized TBA equation acts on the auxiliary variable $z$.
It remains to determine whether the period generating function
itself satisfies a scalar linear $q$-difference equation
of finite order in the modulus.
Writing
\[
\mathscr G(\tau,q):=G(\sqrt{\tau}),
\qquad \tau=T^2,
\]
with the dependence on $q$ made explicit, the question is
whether there exist polynomials $a_j(\tau,q)\in\mathbb C(q)[\tau]$,
not all zero, such that
\[
\sum_{j=0}^{r}
a_j(\tau,q)\,\mathscr G(q^j\tau,q)=0
\]
for some finite $r$, with $q$ treated as an indeterminate.
Our search finds no such equation of order at most eight
with polynomial coefficients of degree at most eighteen.
Determining whether an equation exists beyond these bounds,
or proving its nonexistence, remains an open problem.

\paragraph{Extension to other quantum mirror curves.}
The constant-term definition of quantum periods applies
to much more general Laurent polynomials.
A substantive question is which of these polynomials admit
an analogous correspondence for all Laurent modes,
together with a linearization and continued-fraction
representations.
Natural starting points are nontrivial coefficients in
the present polynomial and other quantum mirror curves.
Such extensions would help identify the algebraic structures
responsible for the correspondence proved here.

\paragraph{Direct computation of the zero mode.}
The algorithms presented in the paper compute the Laurent coefficients
of $\eta(T,z)$ through order $T^{2N}$ using $O(N^2)$
arithmetic operations over $\mathbb C(q)$.
If only the period is required, it is natural to ask whether
the first $N$ nonconstant coefficients of $G(T)$ can be
obtained more efficiently in the same arithmetic model,
without computing all the intermediate modes.
The comparatively simple expression for $1/G(T)$
suggests a possible starting point for this investigation.

 \section*{Acknowledgements}
 We thank L. Chekhov, J. Chuang, L. Gan, V. Golyshev, A. Gorsky, A. Grassi, R. Kashaev, A. Kashani-Poor, M. Kontsevich, A. Mironov, A. Morozov, D.~van Straten, Y. Zenkevich for discussions. The work of Sh. Shakirov and V. Rubtsov is supported by The Ministry of Economic Development of the Russian Federation (IGK 000000C313925P4C0002), agreement No139-15-2025-010.

\bibliographystyle{unsrtnat}
\bibliography{qperbib}

 \newpage
 \appendix

 \section{Differential equation to higher order}
\label{app:ode}

 Here we list the differential equation satisfied by the function $G$ as an expansion around $q=1$ (with $q=e^\epsilon$) to high order, extending \eq{ode3} from the main text. The following equation holds up to $O(\epsilon^{10})$ terms:
 % q=e^\epsilon; all derivatives are with respect to \tau=t^2.
% At epsilon^{2r}, ambiguities are fixed by:
% deg coeff(G)<=r, deg coeff(G')<=r+1,
% and coefficients of G^{(k)}, k>=2, affine in tau.
\begin{equation}
\label{ode9}
\begin{aligned}
&(16\tau^2-\tau)G''+(32\tau-1)G'+4G
\\
&+\epsilon^2\Biggl[
  \frac{\tau}{192}G^{(3)}
  +\frac{1-14\tau}{96}G''
  +\frac{32\tau^2-2\tau-1}{6}G'
  +\frac{32\tau-1}{12}G
\Biggr]
\\
&+\epsilon^4\Biggl[
  \frac{\tau}{2949120}G^{(5)}
  +\frac{1-69\tau}{737280}G^{(4)}
\\
&\hspace{3em}
  +\frac{13(13\tau-1)}{46080}G^{(3)}
  +\frac{1413-3808\tau}{184320}G''
\\
&\hspace{3em}
  -\frac{75776\tau^3-5248\tau^2+312\tau+147}{5760}G'
  -\frac{18944\tau^2-720\tau+45}{2880}G
\Biggr]
\\
&+\epsilon^6\Biggl[
  \frac{\tau}{31708938240}G^{(7)}
  +\frac{1-110\tau}{5284823040}G^{(6)}
\\
&\hspace{3em}
  +\frac{13016\tau-413}{3963617280}G^{(5)}
  +\frac{26189-301664\tau}{1981808640}G^{(4)}
\\
&\hspace{3em}
  +\frac{61864\tau-14471}{30965760}G^{(3)}
  +\frac{44369-36216\tau}{10321920}G''
\\
&\hspace{3em}
  +\frac{61915136\tau^4-5139456\tau^3
         +211264\tau^2-6928\tau-4333}{967680}G'
\\
&\hspace{3em}
  +\frac{15478784\tau^3-801152\tau^2
         +29016\tau-1273}{483840}G
\Biggr]
\\
&+\epsilon^8\Biggl[
  \frac{\tau}{324699527577600}G^{(9)}
  +\frac{2-301\tau}{81174881894400}G^{(8)}
\\
&\hspace{3em}
  +\frac{25992\tau-527}{20293720473600}G^{(7)}
  +\frac{31249-675312\tau}{4058744094720}G^{(6)}
\\
&\hspace{3em}
  +\frac{4154578\tau-397911}{475634073600}G^{(5)}
  +\frac{1054207-5303842\tau}{29727129600}G^{(4)}
\\
&\hspace{3em}
  +\frac{66664016\tau-33018077}{59454259200}G^{(3)}
  +\frac{36550241-11806464\tau}{14863564800}G''
\\
&\hspace{3em}
  -\frac{
    \substack{
      470549528576\tau^5-45423460352\tau^4
      +1702678528\tau^3\\
      {}-47026176\tau^2+2361328\tau+920823
    }
  }{928972800}G'
\\
&\hspace{3em}
  -\frac{
    \substack{
      58818691072\tau^4-3839848448\tau^3\\
      {}+112338688\tau^2-4551904\tau+198915
    }
  }{232243200}G
\Biggr]
=0.
\end{aligned}
\end{equation}

\end{document}